\documentstyle[12pt, leqno, amstex,a4paper]{amsart}

\makeatletter

\begin{document}
\bibliographystyle{amsalpha}

\newcommand\lra{\longrightarrow}
\newcommand\ra{\rightarrow}
\newtheorem{Assumption}{Assumption}[section]
\newtheorem{thm}{Theorem}[section]
\newtheorem{lmm}{Lemma}[section]
\newtheorem{Remark}{Remark}[section]
\newtheorem{Corollary}{Corollary}[section]
\newtheorem{Conjecture}{Conjecture}[section]
\newtheorem{prp}{Proposition}[section]
\newtheorem{exa}{Example}[section]
\newtheorem{dfn}{Definition}[section]
\newtheorem{Problem}{Problem}[section]
\newtheorem{Subsection}{Subsection}[section]
\newtheorem{Proof}{Proof}[section]
\renewcommand{\thesubsection}{\it}

\baselineskip=14pt
\addtocounter{section}{-1}

\title[Coupled Painlev\'e VI systems]{Coupled Painlev\'e VI systems in dimension four with affine Weyl group symmetry\\
 of type $D_6^{(1)}$, II \\}

\author{Yusuke Sasano }
\maketitle

\begin{abstract}
We give a reformulation of a six-parameter family of coupled Pain\-lev\'e VI systems with affine Weyl group symmetry of type  $D_6^{(1)}$ from the viewpoint of its symmetry and holomorphy properties.
\end{abstract}

\section{Introduction}
In \cite{Sasa4,Sasa5}, we proposed a 6-parameter family of four-dimensional coupled Painlev\'e VI systems with affine Weyl group symmetry of type $D_6^{(1)}$. This system can be considered as a genelarization of the Painlev\'e VI system. In this paper, from the viewpoint of its symmetry and holomorphy properties we give a reformulation of this system (cf. \cite{Sasa6}) explicitly given by
\begin{align}\label{1}
\begin{split}
&\frac{dq_1}{dt}=\frac{\partial H}{\partial p_1}, \ \ \frac{dp_1}{dt}=-\frac{\partial H}{\partial q_1}, \ \ \frac{dq_2}{dt}=\frac{\partial H}{\partial p_2}, \ \ \frac{dp_2}{dt}=-\frac{\partial H}{\partial q_2},\\
H =&H_{VI}(q_1,p_1,\eta,t;\alpha_0,\alpha_1,\alpha_2,\alpha_3+2\alpha_4+\alpha_5,\alpha_3+\alpha_6)\\
&+H_{VI}(q_2,p_2,\eta,t;\alpha_0+2\alpha_2+\alpha_3,\alpha_1+\alpha_3,\alpha_4,\alpha_5,\alpha_6)\\
  &+\frac{2(q_1-\eta)q_2\{(q_1-t)p_1+\alpha_2\}\{(q_2-1)p_2+\alpha_4\}}{t(t-1)(t-\eta)} \quad (\eta \in {\Bbb C}-\{0,1\}).
\end{split}
\end{align}
Here $q_1,p_1,q_2,p_2$ denote unknown complex variables, and $\alpha_0,\alpha_1,\dots,\alpha_6$ are complex parameters satisfying the relation $\alpha_0+\alpha_1+2(\alpha_2+\alpha_3+\alpha_4)+\alpha_5+\alpha_6=1$, where the symbol $H_{VI}(q,p,\eta,t;\beta_0,\beta_1,\beta_2,\beta_3,\beta_4)$ is given in Section 2.

If we take the limit $\eta \rightarrow \infty$, we obtain the Hamiltonian system with well-known Hamiltonian $\tilde{H}$ (see \cite{Sasa4})
\begin{align}\label{OCPVI}
\begin{split}
&\frac{dq_1}{dt}=\frac{\partial \tilde{H}}{\partial p_1}, \ \ \frac{dp_1}{dt}=-\frac{\partial \tilde{H}}{\partial q_1}, \ \ \frac{dq_2}{dt}=\frac{\partial \tilde{H}}{\partial p_2}, \ \ \frac{dp_2}{dt}=-\frac{\partial \tilde{H}}{\partial q_2},\\
\tilde{H} &=\tilde{H}_{VI}(q_1,p_1,t;\alpha_0,\alpha_1,\alpha_2,\alpha_3+2\alpha_4+\alpha_5,\alpha_3+\alpha_6)\\
&+\tilde{H}_{VI}(q_2,p_2,t;\alpha_0+\alpha_3,\alpha_1+2\alpha_2+\alpha_3,\alpha_4,\alpha_5,\alpha_6)+\frac{2(q_1-t)p_1q_2\{(q_2-1)p_2+\alpha_4\}}{t(t-1)},
\end{split}
\end{align}
where the symbol $\tilde{H}_{VI}$ is also given in Section 2.

Here we review the holomophy conditions of the system \eqref{OCPVI} (see \cite{Sasa4}). Let us consider a polynomial Hamiltonian system with Hamiltonian $H \in {\Bbb C}(t)[q_1,p_1,q_2,p_2]$. We assume that

$(A1)$ $deg(H)=5$ with respect to $q_1,p_1,q_2,p_2$.

$(A2)$ This system becomes again a polynomial Hamiltonian system in each coordinate system $(x_i,y_i,z_i,w_i) \ (i=0,2,3,4,5,6)${\rm : \rm}
\begin{align}
\begin{split}
&r_0':x_0=-((q_1-t)p_1-\alpha_0)p_1, \ y_0=1/p_1, \ z_0=q_2, \ w_0=p_2, \\
&r_2':x_2=1/q_1, \ y_2=-q_1(q_1p_1+\alpha_2), \ z_2=q_2, \ w_2=p_2, \\
&r_3':x_3=-((q_1-q_2)p_1-\alpha_3)p_1, \ y_3=1/p_1, \ z_3=q_2, \ w_3=p_2+p_1, \\
&r_4':x_4=q_1, \ y_4=p_1, \ z_4=1/q_2, \ w_4=-q_2(q_2p_2+\alpha_4), \\
&r_5':x_5=q_1, \ y_5=p_1, \ z_5=-((q_2-1)p_2-\alpha_5)p_2, \ w_5=1/p_2, \\
&r_6':x_6=q_2, \ y_6=p_1, \ z_6=-p_2(q_2p_2-\alpha_6), \ w_6=1/p_2.
\end{split}
\end{align}

$(A3)$ In addition to the assumption $(A2)$, the Hamiltonian system in the coordinate $r_2$ becomes again a polynomial Hamiltonian system in the coordinate system $(x_1,y_1,z_1,w_1)${\rm : \rm}
\begin{equation}
r_1':x_1=-(x_2y_2-\alpha_1)y_2, \ y_1=1/y_2, \ z_1=z_2, \ w_1=w_2.
\end{equation}
Then such a system coincides with the system \eqref{OCPVI}.

In this paper, we make a reformulation to obtain a clear description of invariant divisors, birational symmetries and holomorphy conditions for the system \eqref{OCPVI}. Our way is stated as follows:
\begin{enumerate}
\item We symmetrize the holomorphy conditions $r_i'$ of the system \eqref{OCPVI}.
\item By using these conditions and {\it polynomiality} of the Hamiltonian, we easily obtain the polynomial Hamiltonian of the system \eqref{1}.
\end{enumerate}

This paper is organized as follows. In Section 2, we give a reformulation of Hamiltonian of $P_{VI}$ and its symmetry and holomorphy. In Section 3, we state our main results for the system of type $D_6^{(1)}$. After we review the notion of accessible singularity in Section 4, we will state the relation between some accessible singularities of the system \eqref{1} and the holomorphy conditions $r_i$ given in Section 3.

\section{Reformulation of $P_{VI}$-case}
The sixth Painlev\'e system can be written as the Hamiltonian system (cf. \cite{Ka,N3})
\begin{align}\label{SPVI}
\begin{split}
&\frac{dq}{dt}=\frac{\partial H_{VI}}{\partial p}, \quad \frac{dp}{dt}=-\frac{\partial H_{VI}}{\partial q},\\
&t(t-1)(t-\eta)H_{VI}(q,p,\eta,t;\alpha_0,\alpha_1,\alpha_2,\alpha_3,\alpha_4)\\
&=q(q-1)(q-\eta)(q-t)p^2+\{\alpha_1(t-\eta)q(q-1)+2\alpha_2q(q-1)(q-\eta)\\
&+\alpha_3(t-1)q(q-\eta)+\alpha_4 t(q-1)(q-\eta)\}p+\alpha_2\{(\alpha_1+\alpha_2)(t-\eta)+\alpha_2(q-1)\\
&+\alpha_3(t-1)+t \alpha_4\}q \quad (\alpha_0+\alpha_1+2\alpha_2+\alpha_3+\alpha_4=1, \quad \eta \in {\Bbb C}-\{0,1\}).
\end{split}
\end{align}
The equation for $q$ is given by
\begin{align}
\begin{split}
\frac{d^2q}{dt^2}=&\frac{1}{2}\left(\frac{1}{q}+\frac{1}{q-1}+\frac{1}{q-t}+\frac{1}{q-\eta}\right)\left(\frac{dq}{dt}\right)^2-\left(\frac{1}{t}+\frac{1}{t-1}+\frac{1}{q-t}+\frac{1}{t-\eta}\right)\frac{dq}{dt}\\
&+\frac{q(q-1)(q-t)(q-\eta)}{t^2(t-1)^2(t-\eta)^2}\{\frac{\alpha_1^2}{2}\frac{\eta(\eta-1)(t-\eta)}{(q-\eta)^2}+\frac{\alpha_4^2}{2}\frac{\eta t}{q^2}\\
&\hspace{4cm}+\frac{\alpha_3^2}{2}\frac{(\eta-1)(1-t)}{(q-1)^2}+\frac{(1-\alpha_0^2)}{2}\frac{t(t-1)(t-\eta)}{(q-t)^2}\}.
\end{split}
\end{align}
If we take the limit $\eta \rightarrow \infty$, we obtain the sixth Painlev\'e system $P_{VI}$ with well-known Hamiltonian:
\begin{align}\label{HVI}
\begin{split}
&\frac{dq}{dt}=\frac{\partial \tilde{H}_{VI}}{\partial p}, \quad \frac{dp}{dt}=-\frac{\partial \tilde{H}_{VI}}{\partial q},\\
&\tilde{H}_{VI}(q,p,t;\delta_0,\delta_1,\delta_2,\delta_3,\delta_4)\\
&=\frac{1}{t(t-1)}[p^2(q-t)(q-1)q-\{(\delta_0-1)(q-1)q+\delta_3(q-t)q\\
&+\delta_4(q-t)(q-1)\}p+\delta_2(\delta_1+\delta_2)q]  \quad (\delta_0+\delta_1+2\delta_2+\delta_3+\delta_4=1),
\end{split}
\end{align}
whose equation for $q$ is given by
\begin{align}
\begin{split}
\frac{d^2q}{dt^2}=&\frac{1}{2}\left(\frac{1}{q}+\frac{1}{q-1}+\frac{1}{q-t}\right)\left(\frac{dq}{dt}\right)^2-\left(\frac{1}{t}+\frac{1}{t-1}+\frac{1}{q-t}\right)\frac{dq}{dt}\\
&+\frac{q(q-1)(q-t)}{t^2(t-1)^2}\{\frac{\alpha_1^2}{2}-\frac{\alpha_4^2}{2}\frac{t}{q^2}-\frac{\alpha_3^2}{2}\frac{(1-t)}{(q-1)^2}+\frac{(1-\alpha_0^2)}{2}\frac{t(t-1)}{(q-t)^2}\}.
\end{split}
\end{align}

The system \eqref{SPVI} has extended affine Weyl group symmetry of type $D_4^{(1)}$, whose generators $s_i,\pi_j$ are given by
\begin{align}\label{D4}
\begin{split}
s_0(q,p,t;\alpha_0,\alpha_1,\dots,\alpha_4) \rightarrow &(q,p-\frac{\alpha_0}{q-t},p,t;-\alpha_0,\alpha_1,\alpha_2+\alpha_0,\alpha_3,\alpha_4),\\
s_1(q,p,t;\alpha_0,\alpha_1,\dots,\alpha_4) \rightarrow &(q,p-\frac{\alpha_1}{q-\eta},t;\alpha_0,-\alpha_1,\alpha_2+\alpha_1,\alpha_3,\alpha_4),\\
s_2(q,p,t;\alpha_0,\alpha_1,\dots,\alpha_4) \rightarrow &(q+\frac{\alpha_2}{p},p,t;\alpha_0+\alpha_2,\alpha_1+\alpha_2,-\alpha_2,\alpha_3+\alpha_2,\alpha_4+\alpha_2),\\
s_3(q,p,t;\alpha_0,\alpha_1,\dots,\alpha_4) \rightarrow &(q,p-\frac{\alpha_3}{q-1},t;\alpha_0,\alpha_1,\alpha_2+\alpha_3,-\alpha_3,\alpha_4),\\
s_4(q,p,t;\alpha_0,\alpha_1,\dots,\alpha_4) \rightarrow &(q,p-\frac{\alpha_4}{q},t;\alpha_0,\alpha_1,\alpha_2+\alpha_4,\alpha_3,-\alpha_4),\\
\pi_1(q,p,t;\alpha_0,\alpha_1,\dots,\alpha_4) \rightarrow &(1-q,-p,1-\eta,1-t;\alpha_0,\alpha_1,\alpha_2,\alpha_4,\alpha_3),\\
\pi_2(q,p,t;\alpha_0,\alpha_1,\dots,\alpha_4) \rightarrow &(\frac{\eta-q}{\eta-1},(1-\eta)p,\frac{\eta}{\eta-1},\frac{\eta-t}{\eta-1};\alpha_0,\alpha_4,\alpha_2,\alpha_3,\alpha_1),
\end{split}
\end{align}

\begin{align*}
\pi_3(q,p,t;\alpha_0,\alpha_1,\dots,\alpha_4) \rightarrow &(\frac{(\eta-1)^2(q-t)}{\{\eta(t-2)+1\}q+(\eta-\eta^2-1)t+\eta^2},\\
&(1-t)p+\frac{(q-1)\{(q-1)p+\alpha_2\}\{\eta(t-2)+1\}}{(\eta-1)^2(t-1)}\\
&+\frac{(q-t)\{(q-t)p+\alpha_2\}\{\eta(t-2)+1\}}{\eta(t-1)(t-\eta)},\\
&1-\eta,\frac{(\eta-1)^2t}{t-\eta t+\eta^2(t-1)};\alpha_4,\alpha_1,\alpha_2,\alpha_3,\alpha_0).
\end{align*}

Let us consider a polynomial Hamiltonian system with Hamiltonian $H \in {\Bbb C}(t)[q,p]$. We assume that

$(A1)$ $deg(H)=6$ with respect to $q,p$.

$(A2)$ This system becomes again a polynomial Hamiltonian system in each coordinate $r_i \ (j=0,1,\dots,4)$:
\begin{align}\label{holoPVI}
\begin{split}
&r_0:x_0=-((q-t)p-\alpha_0)p,\ y_0=\frac{1}{p},\\
&r_1:x_1=-((q-\eta)p-\alpha_1)p, \ y_1=\frac{1}{p}, \\
&r_2:x_2=\frac{1}{q}, \ y_2=-(qp+\alpha_2)q,\\
&r_3:x_3=-((q-1)p-\alpha_3)p,\ y_3=\frac{1}{p},\\
&r_4:x_4=-(qp-\alpha_4)p,\ y_4=\frac{1}{p}.
\end{split}
\end{align}
Then such a system coincides with the system \eqref{SPVI}.

The phase space of the system \eqref{SPVI} (resp. \eqref{HVI}) can be characterized by the rational surface of type $D_4^{(1)}$ (see \cite{Oka,Saito,Sakai}). The below figure denotes the accessible singular points and the resolution process for each system.

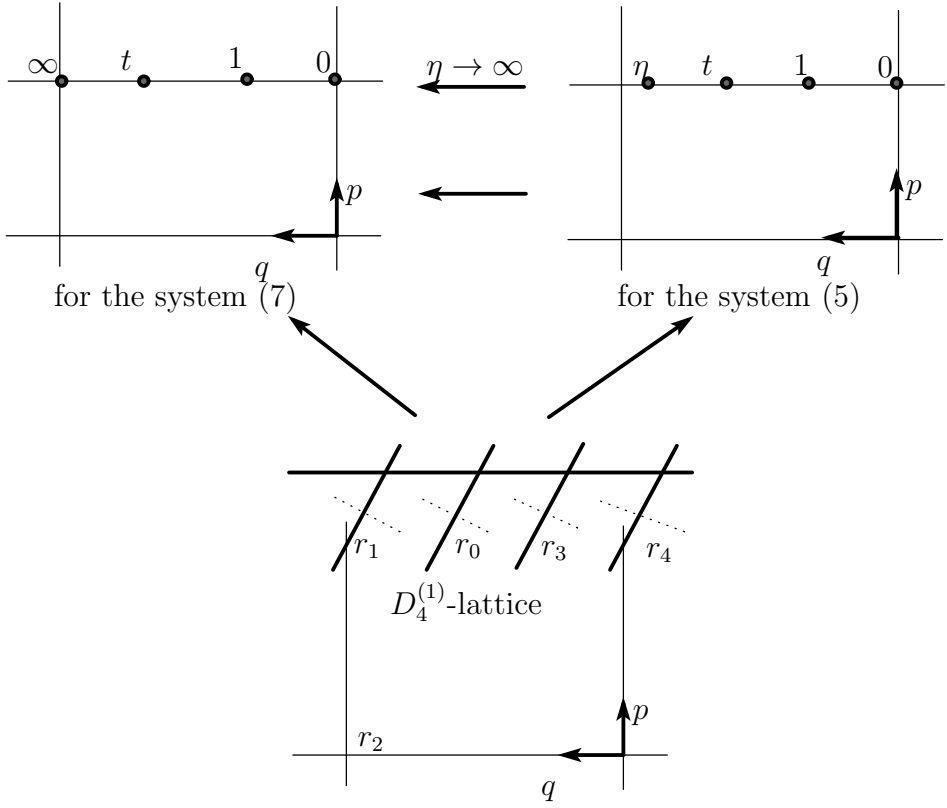
\begin{figure}[h]
\unitlength 0.1in
\begin{picture}(49.10,41.20)(6.10,-45.20)
%
\special{pn 8}%
\special{pa 620 810}%
\special{pa 2580 810}%
\special{fp}%
%
\special{pn 8}%
\special{pa 890 400}%
\special{pa 890 1780}%
\special{fp}%
%
\special{pn 8}%
\special{pa 610 1620}%
\special{pa 2570 1620}%
\special{fp}%
%
\special{pn 8}%
\special{pa 2340 420}%
\special{pa 2340 1800}%
\special{fp}%
%
\special{pn 20}%
\special{sh 0.600}%
\special{ar 900 810 25 25  0.0000000 6.2831853}%
%
\special{pn 20}%
\special{sh 0.600}%
\special{ar 1330 810 25 25  0.0000000 6.2831853}%
%
\special{pn 20}%
\special{sh 0.600}%
\special{ar 1870 800 25 25  0.0000000 6.2831853}%
%
\special{pn 20}%
\special{sh 0.600}%
\special{ar 2330 800 25 25  0.0000000 6.2831853}%
%
\special{pn 8}%
\special{pa 3560 830}%
\special{pa 5520 830}%
\special{fp}%
%
\special{pn 8}%
\special{pa 3830 420}%
\special{pa 3830 1800}%
\special{fp}%
%
\special{pn 8}%
\special{pa 3550 1640}%
\special{pa 5510 1640}%
\special{fp}%
%
\special{pn 8}%
\special{pa 5280 440}%
\special{pa 5280 1820}%
\special{fp}%
%
\special{pn 20}%
\special{sh 0.600}%
\special{ar 3970 820 25 25  0.0000000 6.2831853}%
%
\special{pn 20}%
\special{sh 0.600}%
\special{ar 4380 820 25 25  0.0000000 6.2831853}%
%
\special{pn 20}%
\special{sh 0.600}%
\special{ar 4810 820 25 25  0.0000000 6.2831853}%
%
\special{pn 20}%
\special{sh 0.600}%
\special{ar 5270 820 25 25  0.0000000 6.2831853}%
%
\special{pn 20}%
\special{pa 2340 1620}%
\special{pa 2340 1350}%
\special{fp}%
\special{sh 1}%
\special{pa 2340 1350}%
\special{pa 2320 1417}%
\special{pa 2340 1403}%
\special{pa 2360 1417}%
\special{pa 2340 1350}%
\special{fp}%
%
\special{pn 20}%
\special{pa 2340 1620}%
\special{pa 2030 1620}%
\special{fp}%
\special{sh 1}%
\special{pa 2030 1620}%
\special{pa 2097 1640}%
\special{pa 2083 1620}%
\special{pa 2097 1600}%
\special{pa 2030 1620}%
\special{fp}%
\put(19.1000,-18.3000){\makebox(0,0)[lb]{$q$}}%
\put(23.9000,-14.3000){\makebox(0,0)[lb]{$p$}}%
\put(8.6000,-19.9000){\makebox(0,0)[lb]{for the system \eqref{HVI}}}%
\put(22.3000,-7.6000){\makebox(0,0)[lb]{0}}%
\put(17.7000,-7.4000){\makebox(0,0)[lb]{1}}%
\put(12.1000,-7.5000){\makebox(0,0)[lb]{$t$}}%
\put(7.2000,-7.7000){\makebox(0,0)[lb]{$\infty$}}%
\put(38.9000,-7.8000){\makebox(0,0)[lb]{$\eta$}}%
%
\special{pn 20}%
\special{pa 3320 840}%
\special{pa 2790 840}%
\special{fp}%
\special{sh 1}%
\special{pa 2790 840}%
\special{pa 2857 860}%
\special{pa 2843 840}%
\special{pa 2857 820}%
\special{pa 2790 840}%
\special{fp}%
%
\special{pn 20}%
\special{pa 3330 1400}%
\special{pa 2800 1400}%
\special{fp}%
\special{sh 1}%
\special{pa 2800 1400}%
\special{pa 2867 1420}%
\special{pa 2853 1400}%
\special{pa 2867 1380}%
\special{pa 2800 1400}%
\special{fp}%
\put(28.1000,-7.8000){\makebox(0,0)[lb]{$\eta \rightarrow \infty$}}%
\put(38.1000,-19.9000){\makebox(0,0)[lb]{for the system \eqref{SPVI}}}%
%
\special{pn 20}%
\special{pa 5280 1630}%
\special{pa 4910 1630}%
\special{fp}%
\special{sh 1}%
\special{pa 4910 1630}%
\special{pa 4977 1650}%
\special{pa 4963 1630}%
\special{pa 4977 1610}%
\special{pa 4910 1630}%
\special{fp}%
%
\special{pn 20}%
\special{pa 5270 1630}%
\special{pa 5270 1300}%
\special{fp}%
\special{sh 1}%
\special{pa 5270 1300}%
\special{pa 5250 1367}%
\special{pa 5270 1353}%
\special{pa 5290 1367}%
\special{pa 5270 1300}%
\special{fp}%
\put(48.5000,-18.0000){\makebox(0,0)[lb]{$q$}}%
\put(53.2000,-14.0000){\makebox(0,0)[lb]{$p$}}%
\put(51.7000,-7.9000){\makebox(0,0)[lb]{0}}%
\put(47.3000,-7.8000){\makebox(0,0)[lb]{1}}%
\put(42.5000,-7.8000){\makebox(0,0)[lb]{$t$}}%
%
\special{pn 8}%
\special{pa 2390 3120}%
\special{pa 2390 4500}%
\special{fp}%
%
\special{pn 8}%
\special{pa 2110 4340}%
\special{pa 4070 4340}%
\special{fp}%
%
\special{pn 8}%
\special{pa 3840 3140}%
\special{pa 3840 4520}%
\special{fp}%
%
\special{pn 20}%
\special{pa 3840 4340}%
\special{pa 3840 4070}%
\special{fp}%
\special{sh 1}%
\special{pa 3840 4070}%
\special{pa 3820 4137}%
\special{pa 3840 4123}%
\special{pa 3860 4137}%
\special{pa 3840 4070}%
\special{fp}%
%
\special{pn 20}%
\special{pa 3840 4340}%
\special{pa 3530 4340}%
\special{fp}%
\special{sh 1}%
\special{pa 3530 4340}%
\special{pa 3597 4360}%
\special{pa 3583 4340}%
\special{pa 3597 4320}%
\special{pa 3530 4340}%
\special{fp}%
\put(34.1000,-45.5000){\makebox(0,0)[lb]{$q$}}%
\put(38.9000,-41.5000){\makebox(0,0)[lb]{$p$}}%
%
\special{pn 20}%
\special{pa 2320 3370}%
\special{pa 2670 2720}%
\special{fp}%
%
\special{pn 20}%
\special{pa 2810 3370}%
\special{pa 3160 2720}%
\special{fp}%
%
\special{pn 20}%
\special{pa 3280 3360}%
\special{pa 3630 2710}%
\special{fp}%
%
\special{pn 20}%
\special{pa 3770 3370}%
\special{pa 4120 2720}%
\special{fp}%
%
\special{pn 20}%
\special{pa 2090 2860}%
\special{pa 4200 2860}%
\special{fp}%
%
\special{pn 8}%
\special{pa 2330 2990}%
\special{pa 2670 3170}%
\special{dt 0.045}%
\special{pa 2670 3170}%
\special{pa 2669 3170}%
\special{dt 0.045}%
%
\special{pn 8}%
\special{pa 2790 3000}%
\special{pa 3130 3170}%
\special{dt 0.045}%
\special{pa 3130 3170}%
\special{pa 3129 3170}%
\special{dt 0.045}%
%
\special{pn 8}%
\special{pa 3270 3000}%
\special{pa 3600 3150}%
\special{dt 0.045}%
\special{pa 3600 3150}%
\special{pa 3599 3150}%
\special{dt 0.045}%
%
\special{pn 8}%
\special{pa 3720 3010}%
\special{pa 4160 3170}%
\special{dt 0.045}%
\special{pa 4160 3170}%
\special{pa 4159 3170}%
\special{dt 0.045}%
%
\special{pn 20}%
\special{pa 2750 2560}%
\special{pa 2110 2060}%
\special{fp}%
\special{sh 1}%
\special{pa 2110 2060}%
\special{pa 2150 2117}%
\special{pa 2152 2093}%
\special{pa 2175 2085}%
\special{pa 2110 2060}%
\special{fp}%
%
\special{pn 20}%
\special{pa 3450 2570}%
\special{pa 4180 2060}%
\special{fp}%
\special{sh 1}%
\special{pa 4180 2060}%
\special{pa 4114 2082}%
\special{pa 4136 2091}%
\special{pa 4137 2115}%
\special{pa 4180 2060}%
\special{fp}%
%
\put(42.0000,-23.4000){\makebox(0,0)[lb]{}}%
\put(26.2000,-36.0000){\makebox(0,0)[lb]{$D_4^{(1)}$-lattice}}%
\put(24.5000,-42.9000){\makebox(0,0)[lb]{$r_2$}}%
\put(39.6000,-33.0000){\makebox(0,0)[lb]{$r_4$}}%
\put(34.1000,-32.9000){\makebox(0,0)[lb]{$r_3$}}%
\put(29.6000,-32.8000){\makebox(0,0)[lb]{$r_0$}}%
\put(24.3000,-32.8000){\makebox(0,0)[lb]{$r_1$}}%
\end{picture}%
\label{HPVI}
\caption{Each figure denotes the Hirzebruch surface. Each bullet denotes the accessible singular point of each system. It is well-known that each point can be resolved by blowing-up at two times (see \cite{Oka,Saito,Sakai}). By these transformations, we obtain the rational surface of type $D_4^{(1)}$ for each system.}
\end{figure}

We remark that the system \eqref{SPVI} has the following invariant divisors\rm{:\rm}
\begin{center}
\begin{tabular}{|c|c|} \hline
parameter's relation & invariant divisors \\ \hline
$\alpha_0=0$ & $f_0:=q-t$  \\ \hline
$\alpha_1=0$ & $f_1:=q-\eta$  \\ \hline
$\alpha_2=0$ & $f_2:=p$ \\ \hline
$\alpha_3=0$ & $f_3:=q-1$  \\ \hline
$\alpha_4=0$ & $f_4:=q$   \\ \hline
\end{tabular}
\end{center}

\section{The case of type $D_6^{(1)}$}

\begin{thm}\label{1.1}
The system \eqref{1} admits extended affine Weyl group symmetry of type $D_6^{(1)}$ as the group of its B{\"a}cklund transformations, whose generators $s_i,{\pi}_j$ are explicitly given as follows{\rm : \rm}with the notation $(*):=(q_1,p_1,q_2,p_2,\eta,t;\alpha_0,\alpha_1,\dots,\alpha_6),$
\begin{align*}
        s_0: &(*) \rightarrow (q_1,p_1-\frac{\alpha_0}{q_1-t},q_2,p_2,\eta,t;-\alpha_0,\alpha_1,\alpha_2+\alpha_0,\alpha_3,\alpha_4,\alpha_5,\alpha_6),\\
        s_1: &(*) \rightarrow (q_1,p_1-\frac{\alpha_1}{q_1-\eta},q_2,p_2,\eta,t;\alpha_0,-\alpha_1,\alpha_2+\alpha_1,\alpha_3,\alpha_4,\alpha_5,\alpha_6),\\
        s_2: &(*) \rightarrow (q_1+\frac{\alpha_2}{p_1},p_1,q_2,p_2,\eta,t;\alpha_0+\alpha_2,\alpha_1+\alpha_2,-\alpha_2,\alpha_3+\alpha_2,\alpha_4,\alpha_5,\alpha_6), \\
        s_3: &(*) \rightarrow (q_1,p_1-\frac{\alpha_3}{q_1-q_2},q_2,p_2+\frac{\alpha_3}{q_1-q_2},\eta,t;\alpha_0,\alpha_1,\alpha_2+\alpha_3,-\alpha_3,\alpha_4+\alpha_3,\alpha_5,\alpha_6), \\
        s_4: &(*) \rightarrow (q_1,p_1,q_2+\frac{\alpha_4}{p_2},p_2,\eta,t;\alpha_0,\alpha_1,\alpha_2,\alpha_3+\alpha_4,-\alpha_4,\alpha_5+\alpha_4,\alpha_6+\alpha_4),\\
        s_5: &(*) \rightarrow (q_1,p_1,q_2,p_2-\frac{\alpha_5}{q_2-1},\eta,t;\alpha_0,\alpha_1,\alpha_2,\alpha_3,\alpha_4+\alpha_5,-\alpha_5,\alpha_6), \\
        s_6: &(*) \rightarrow (q_1,p_1,q_2,p_2-\frac{\alpha_6}{q_2},\eta,t;\alpha_0,\alpha_1,\alpha_2,\alpha_3,\alpha_4+\alpha_6,\alpha_5,-\alpha_6),
        \end{align*}
        \begin{align*}
        {\pi}_1: &(*) \rightarrow\\
&(\frac{(t-1)q_1}{t-q_1-\eta t+\eta tq_1},\frac{(-t+q_1+\eta t-\eta tq_1)(tp_1-q_1p_1-\alpha_2-\eta tp_1+\eta tq_1p_1+\alpha_2 \eta t)}{t(t-1)(\eta-1)},\\
        &\frac{(t-1)q_2}{t-q_2-\eta t+\eta tq_2},\frac{(-t+q_2+\eta t-\eta tq_2)(tp_2-q_2p_2-\alpha_4-\eta tp_2+\eta tq_2p_2+\alpha_4 \eta t)}{t(t-1)(\eta-1)},\\
        &\frac{1}{\eta},\frac{\eta(t-1)}{t-\eta-\eta t+\eta^2 t};\alpha_1,\alpha_0,\alpha_2,\alpha_3,\alpha_4,\alpha_5,\alpha_6),\\
        {\pi}_2: &(*) \rightarrow  (1-q_1,-p_1,1-q_2,-p_2,1-\eta,1-t;\alpha_0,\alpha_1,\alpha_2,\alpha_3,\alpha_4,\alpha_6,\alpha_5), \\
        {\pi}_3: &(*) \rightarrow\\
        &(\frac{t(q_2-\eta)}{t(q_2-\eta)+\eta^2 (t-q_2)},\\
        &\frac{(t(q_2-\eta)+\eta^2 (t-q_2))(t(q_2-\eta)p_2+\alpha_4(t-\eta^2)+\eta^2(t-q_2)p_2)}{t\eta^2(t-\eta)},\\
&\frac{t(q_1-\eta)}{t(q_1-\eta)+\eta^2 (t-q_1)},\\
&\frac{(t(q_1-\eta)+\eta^2 (t-q_1))(t(q_1-\eta)p_1+\alpha_2(t-\eta^2)+\eta^2(t-q_1)p_1)}{t\eta^2(t-\eta)},\\
        &-\frac{1}{\eta-1},-\frac{(\eta-1)t}{t-\eta t+\eta^2(t-1)};\alpha_5,\alpha_6,\alpha_4,\alpha_3,\alpha_2,\alpha_0,\alpha_1).
\end{align*}
\end{thm}
We note that these transformations $s_i,\pi_j$ are birational and symplectic.

\begin{figure}
\unitlength 0.1in
\begin{picture}(41.56,10.30)(22.70,-17.21)
%
\special{pn 20}%
\special{ar 2772 1010 236 123  1.5707963 6.2831853}%
\special{ar 2772 1010 236 123  0.0000000 1.5284488}%
%
\special{pn 20}%
\special{ar 2772 1597 236 124  1.5707963 6.2831853}%
\special{ar 2772 1597 236 124  0.0000000 1.5284488}%
%
\special{pn 20}%
\special{ar 3548 1313 236 125  1.5707963 6.2831853}%
\special{ar 3548 1313 236 125  0.0000000 1.5284488}%
%
\special{pn 20}%
\special{ar 5316 1313 235 125  1.5707963 6.2831853}%
\special{ar 5316 1313 235 125  0.0000000 1.5282688}%
%
\special{pn 20}%
\special{ar 6190 1004 236 124  1.5707963 6.2831853}%
\special{ar 6190 1004 236 124  0.0000000 1.5326792}%
%
\special{pn 20}%
\special{ar 6190 1597 236 124  1.5707963 6.2831853}%
\special{ar 6190 1597 236 124  0.0000000 1.5326792}%
%
\special{pn 20}%
\special{pa 2978 1061}%
\special{pa 3362 1226}%
\special{fp}%
%
\special{pn 20}%
\special{pa 2988 1551}%
\special{pa 3381 1401}%
\special{fp}%
%
\special{pn 20}%
\special{ar 4452 1318 235 124  1.5707963 6.2831853}%
\special{ar 4452 1318 235 124  0.0000000 1.5282688}%
%
\special{pn 20}%
\special{pa 3784 1318}%
\special{pa 4196 1318}%
\special{fp}%
%
\special{pn 20}%
\special{pa 4687 1318}%
\special{pa 5061 1318}%
\special{fp}%
%
\special{pn 20}%
\special{pa 5523 1236}%
\special{pa 5983 1081}%
\special{fp}%
%
\special{pn 20}%
\special{pa 5493 1391}%
\special{pa 5964 1556}%
\special{fp}%
%
\put(22.7000,-9.9200){\makebox(0,0)[lb]{}}%
\put(25.8500,-16.6800){\makebox(0,0)[lb]{$q_1-t$}}%
\put(34.3000,-13.5800){\makebox(0,0)[lb]{$p_1$}}%
\put(42.5000,-13.7000){\makebox(0,0)[lb]{$q_1-q_2$}}%
\put(51.9900,-13.8300){\makebox(0,0)[lb]{$p_2$}}%
\put(60.3600,-16.5100){\makebox(0,0)[lb]{$q_2-1$}}%
\put(60.8600,-10.6800){\makebox(0,0)[lb]{$q_2$}}%
\put(25.8500,-10.7300){\makebox(0,0)[lb]{$q_1-\eta$}}%
\put(26.0500,-14.6600){\makebox(0,0)[lb]{$\alpha_0$}}%
\put(26.0500,-8.6700){\makebox(0,0)[lb]{$\alpha_1$}}%
\put(34.1000,-11.5400){\makebox(0,0)[lb]{$\alpha_2$}}%
\put(51.3900,-11.6700){\makebox(0,0)[lb]{$\alpha_4$}}%
\put(60.0300,-14.6000){\makebox(0,0)[lb]{$\alpha_5$}}%
\put(60.2300,-8.6100){\makebox(0,0)[lb]{$\alpha_6$}}%
\put(41.9600,-11.5800){\makebox(0,0)[lb]{$\alpha_3$}}%
\end{picture}%
\label{D65}
\caption{This figure denotes Dynkin diagram of type $D_6^{(1)}$.}
\end{figure}
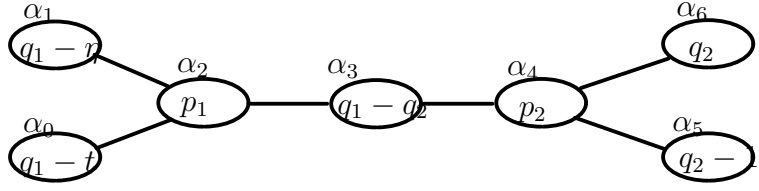

\begin{thm}\label{1.2}
Let us consider a polynomial Hamiltonian system with Hamiltonian $H \in {\Bbb C}(t)[q_1,p_1,q_2,p_2]$. We assume that

$(A1)$ $deg(H)=6$ with respect to $q_1,p_1,q_2,p_2$.

$(A2)$ This system becomes again a polynomial Hamiltonian system in each coordinate system $(x_i,y_i,z_i,w_i) \ (i=0,1,\dots ,6)${\rm : \rm}
\begin{align}
\begin{split}
&r_0:x_0=-((q_1-t)p_1-\alpha_0)p_1, \ y_0=1/p_1, \ z_0=q_2, \ w_0=p_2, \\
&r_1:x_1=-((q_1-\eta)p_1-\alpha_1)p_1, \ y_1=1/p_1, \ z_1=q_2, \ w_1=p_2 \quad (\eta \in {\Bbb C}-\{0,1\}), \\
&r_2:x_2=1/q_1, \ y_2=-q_1(q_1p_1+\alpha_2), \ z_2=q_2, \ w_2=p_2, \\
&r_3:x_3=-((q_1-q_2)p_1-\alpha_3)p_1, \ y_3=1/p_1, \ z_3=q_2, \ w_3=p_2+p_1, \\
&r_4:x_4=q_1, \ y_4=p_1, \ z_4=1/q_2, \ w_4=-q_2(q_2p_2+\alpha_4), \\
&r_5:x_5=q_1, \ y_5=p_1, \ z_5=-((q_2-1)p_2-\alpha_5)p_2, \ w_5=1/p_2, \\
&r_6:x_6=q_1, \ y_6=p_1, \ z_6=-p_2(q_2p_2-\alpha_6), \ w_6=1/p_2.
\end{split}
\end{align}
Then such a system coincides with the system \eqref{1}.
\end{thm}

\begin{prp}
The system \eqref{1} has the following invariant divisors\rm{:\rm}
\begin{center}
\begin{tabular}{|c|c|} \hline
parameter's relation & invariant divisors \\ \hline
$\alpha_0=0$ & $f_0:=q_1-t$  \\ \hline
$\alpha_1=0$ & $f_1:=q_1-\eta$  \\ \hline
$\alpha_2=0$ & $f_2:=p_1$  \\ \hline
$\alpha_3=0$ & $f_3:=q_1-q_2$  \\ \hline
$\alpha_4=0$ & $f_4:=p_2$   \\ \hline
$\alpha_5=0$ & $f_5:=q_2-1$   \\ \hline
$\alpha_6=0$ & $f_6:=q_2$   \\ \hline
\end{tabular}
\end{center}
\end{prp}

\section{Accessible singularities }

Let us review the notion of accessible singularity. Let $B$ be a connected open domain in $\Bbb C$ and $\pi : {\mathcal W} \longrightarrow B$ a smooth proper holomorphic map. We assume that ${\mathcal H} \subset {\mathcal W}$ is a normal crossing divisor which is flat over $B$. Let us consider a rational vector field $\tilde v$ on $\mathcal W$ satisfying the condition
\begin{equation*}
\tilde v \in H^0({\mathcal W},\Theta_{\mathcal W}(-\log{\mathcal H})({\mathcal H})).
\end{equation*}
Fixing $t_0 \in B$ and $P \in {\mathcal W}_{t_0}$, we can take a local coordinate system $(x_1,x_2,\dots ,x_n)$ of ${\mathcal W}_{t_0}$ centered at $P$ such that ${\mathcal H}_{\rm smooth \rm}$ can be defined by the local equation $x_1=0$.
Since $\tilde v \in H^0({\mathcal W},\Theta_{\mathcal W}(-\log{\mathcal H})({\mathcal H}))$, we can write down the vector field $\tilde v$ near $P=(0,0,\dots ,0,t_0)$ as follows:
\begin{equation}
\tilde v= \frac{\partial}{\partial t}+a_1 
\frac{\partial}{\partial x_1}+\frac{a_2}{x_1} 
\frac{\partial}{\partial x_2}+.....+\frac{a_n}{x_1} 
\frac{\partial}{\partial x_n}.
\end{equation}
This vector field defines the following system of differential equations
\begin{equation}\label{39}
  \left\{
  \begin{aligned}
   \frac{dx_1}{dt} &=a_1(x_1,x_2,....,x_n,t),\\
   \frac{dx_2}{dt} &=\frac{a_2(x_1,x_2,....,x_n,t)}{x_1},\\
   .\\
   .\\
   .\\
   \frac{dx_n}{dt} &=\frac{a_n(x_1,x_2,....,x_n,t)}{x_1}.
   \end{aligned}
  \right. 
\end{equation}
Here $a_i(x_1,x_2,....,x_n,t), \ i=1,2,\dots ,n,$ are holomorphic functions defined near $P=(0,\dots ,0,t_0).$

\begin{dfn}\label{Def}
With the above notation, assume that the rational vector field $\tilde v$ on $\mathcal W$ satisfies the condition
$$
(A) \quad \tilde v \in H^0({\mathcal W},\Theta_{\mathcal W}(-\log{\mathcal H})({\mathcal H})).
$$
We say that $\tilde v$ has an {\it accessible singularity} at $P=(0,0,\dots ,0,t_0)$ if
$$
x_1=0 \ {\rm and \rm} \ a_i(0,0,....,0,t_0)=0 \ {\rm for \rm} \ {\rm every \rm} \ i, \ 2 \leq i \leq n.
$$
\end{dfn}

If $P \in {\mathcal H}_{{\rm smooth \rm}}$ is not an accessible singularity, all solutions of the ordinary differential equation passing through $P$ are vertical solutions, that is, the solutions are contained in the fiber ${\mathcal W}_{t_0}$ over $t=t_0$. If $P \in {\mathcal H}_{\rm smooth \rm}$ is an accessible singularity, there may be a solution of \eqref{39} which passes through $P$ and goes into the interior ${\mathcal W}-{\mathcal H}$ of ${\mathcal W}$.

Here we review the notion of {\it local index}. Let $v$ be an algebraic vector field with an accessible singular point $\overrightarrow{p}=(0,0,\ldots,0)$ and $(x_1,x_2,\ldots,x_n)$ be a coordinate system in a neighborhood centered at $\overrightarrow{p}$. Assume that the system associated with $v$ near $\overrightarrow{p}$ can be written as
\begin{align}\label{b}
\begin{split}
\frac{d}{dt}Q\begin{pmatrix}
             x_1 \\
             x_2 \\
             \vdots\\
             x_n
             \end{pmatrix}=\frac{1}{x_1}\{Q\begin{bmatrix}
             a_1 & & & &  \\
             & a_2 & & & \\
             & & \ddots & & \\
             & & & & a_n
             \end{bmatrix}Q^{-1}{\cdot}Q\begin{pmatrix}
             x_1 \\
             x_2 \\
             \vdots\\
             x_n
             \end{pmatrix}+\begin{pmatrix}
             x_1f_1(x_1,x_2,\ldots,x_n,t) \\
             f_2(x_1,x_2,\ldots,x_n,t) \\
             \vdots\\
             f_n(x_1,x_2,\ldots,x_n,t)
             \end{pmatrix}\},\\
              (f_i \in {\Bbb C}(t)[x_1,\ldots,x_n], \ Q \in GL(n,{\Bbb C}(t)),a_i \in {\Bbb C}(t))
             \end{split}
             \end{align}
where $f_1(x_1,x_2,\ldots,x_n,t)$ is a polynomial which vanishes at $\overrightarrow{p}$ and $f_i(x_1,x_2,\ldots,x_n,t)$, $i=2,3,\ldots,n$ are polynomials of order at least 2 in $x_1,x_2,\ldots,x_n$. We call ordered set of the eigenvalues $(a_1,a_2,\ldots,a_n)$ {\it local index} at $\overrightarrow{p}$ if
\begin{equation}
(1,a_2/a_1,\ldots,a_n/a_1) \in {\Bbb Z}^n.
\end{equation}
We remark that if each component of $(1,a_2/a_1,\ldots,a_n/a_1)$ has the same sign, we may resolve the accessible singularity by blowing-up finitely many times. However, when different signs appear, we may need to both blow up and blow down.

\begin{exa}
For the Noumi-Yamada system of type $A_4^{(1)}$, its local index can be defined at each accessible singular point (cf. \cite{Ta}).
\end{exa}

\section{On some Hamiltonian structures of the system \eqref{1}}
In this section, we will give the holomorphy conditions $r_i \ (i=0,1,\ldots,6)$ by resolving some accessible singular loci of the system \eqref{1}. Each of them contains a 3-parameter family of meromorphic solutions.

In order to consider the singularity analysis for the system \eqref{1}, as a compactification of ${\Bbb C}^4$ which is the phase space of the system \eqref{1}, at first we take a 4-dimensional projective space ${\Bbb P}^4$. In this space the rational vector field $\tilde v$ associated with the system \eqref{1} satisfies the condition:
$$
\tilde v \in H^0({{\Bbb P}^4},\Theta_{{\Bbb P}^4}(-\log{H})(3H)),
$$
where $H$ denotes the boundary divisor $H \cong {\Bbb P}^3$.
To calculate its accessible singularities, we must replace the compactification of ${\Bbb C}^4$ with the condition $(A)$ given in Section 4.
In this paper, we present complex manifold $\mathcal S$ obtained by gluing twelve copies $U_j \cong {\Bbb C}^4 \ni (X_j,Y_j,Z_j,W_j), \ j=0,1,\ldots,11$:
\begin{equation*}
U_j \times B={\Bbb C}^4 \times B \ni (X_j,Y_j,Z_j,W_j,t) \ (j=0,1,\ldots,11)
\end{equation*}
via the following birational transformations:
\begin{align}
\begin{split}
0)&X_0=q_1, \quad Y_0=p_1, \quad Z_0=q_2, \quad W_0=p_2,\\
1)&X_1=1/q_1, \quad Y_1=-(q_1p_1+\alpha_2)q_1, \quad Z_1=q_2, \quad W_1=p_2,\\
2)&X_2=q_1, \quad Y_2=p_1, \quad Z_2=1/q_2, \quad W_2=-(q_2p_2+\alpha_4)q_2,\\
3)&X_3=q_1, \quad Y_3=1/p_1, \quad Z_3=q_2, \quad W_3=p_2/p_1,\\
4)&X_4=q_1, \quad Y_4=p_1/p_2, \quad Z_4=q_2, \quad W_4=1/p_2,\\
5)&X_5=1/q_1, \quad Y_5=-(q_1p_1+\alpha_2)q_1, \quad Z_5=1/q_2, \quad W_5=-(q_2p_2+\alpha_4)q_2,\\
6)&X_6=1/q_1, \quad Y_6=-\frac{1}{(q_1p_1+\alpha_2)q_1}, \quad Z_6=q_2, \quad W_6=-\frac{p_2}{(q_1p_1+\alpha_2)q_1},\\
7)&X_7=1/q_1, \quad Y_7=-\frac{(q_1p_1+\alpha_2)q_1}{p_2}, \quad Z_7=q_2, \quad W_7=1/p_2,\\
8)&X_8=1/q_1, \quad Y_8=-\frac{1}{(q_1p_1+\alpha_2)q_1}, \quad Z_8=1/q_2, \quad W_8=\frac{(q_2p_2+\alpha_4)q_2}{(q_1p_1+\alpha_2)q_1},\\
9)&X_9=1/q_1, \quad Y_9=\frac{(q_1p_1+\alpha_2)q_1}{(q_2p_2+\alpha_4)q_2}, \quad Z_9=1/q_2, \quad W_9=-\frac{1}{(q_2p_2+\alpha_4)q_2},\\
10)&X_{10}=q_1, \quad Y_{10}=1/p_1, \quad Z_{10}=1/q_2, \quad W_{10}=-\frac{(q_2p_2+\alpha_4)q_2}{p_1},\\
11)&X_{11}=q_1, \quad Y_{11}=-\frac{p_1}{(q_2p_2+\alpha_4)q_2}, \quad Z_{11}=1/q_2, \quad W_{11}=-\frac{1}{(q_2p_2+\alpha_4)q_2}.
\end{split}
\end{align}
The restriction $\{(q_1,p_1,q_2,p_2)|q_2=p_2=0\}$ (resp. $\{(q_1,p_1,q_2,p_2)|q_1=p_1=0\}$) of this manifold $\mathcal S$ is a Hirzebruch surface, respectively. Thus, it can be considered as a generalization of Hirzebruch surface. We remark that this generalization of the Hirzebruch surface is different from the one given by H. Kimura (see \cite{K}).

The canonical divisor $K_{\mathcal S}$ of $\mathcal S$ is given by
\begin{align}
\begin{split}
K_{\mathcal S}&=-3{\mathcal H}\\
&=\bigcup_{i \in \{3,6,8,10\}}\{(X_i,Y_i,Z_i,W_i) \in U_i|Y_i=0\} \cup \bigcup_{j \in \{4,7,9,11\}}\{(X_j,Y_j,Z_j,W_j) \in U_j|W_j=0\},
\end{split}
\end{align}
and satisfies the following relations:
\begin{equation}
  \left\{
  \begin{aligned}
   dX_j \wedge dY_j \wedge dZ_j \wedge dW_j &=dq_1 \wedge dp_1 \wedge dq_2 \wedge dp_2 \ (j=1,2,5),\\
   dX_3 \wedge dY_3 \wedge dZ_3 \wedge dW_3 &=-\frac{1}{p_1^3}dq_1 \wedge dp_1 \wedge dq_2 \wedge dp_2,\\
   dX_6 \wedge dY_6 \wedge dZ_6 \wedge dW_6 &=-\frac{1}{Y_1^3}dX_1 \wedge dY_1 \wedge dZ_1 \wedge dW_1,\\
   dX_8 \wedge dY_8 \wedge dZ_8 \wedge dW_8 &=-\frac{1}{Y_5^3}dX_5 \wedge dY_5 \wedge dZ_5 \wedge dW_5.\\
   \end{aligned}
  \right. 
\end{equation}
We note that the transformation
\begin{equation}
\pi:(q_1,p_1,q_2,p_2;\alpha_2,\alpha_4) \rightarrow (q_2,p_2,q_1,p_1;\alpha_4,\alpha_2)
\end{equation}
is an automorphism of $\mathcal S$.

It is easy to see that each patching data $(X_i,Y_i,Z_i,W_i) \ (i=1,2,5)$ is birational and symplectic, moreover the system \eqref{1} becomes again a {\it polynomial} Hamiltonian system in each coordinate system.

\begin{prp}
After a series of explicit blowing-ups and blowing-downs of ${\Bbb P}^4$, we obtain the smooth projective 4-fold $\mathcal S$ and a birational morphism $\varphi:{\mathcal S} \dots \rightarrow {\Bbb P}^4$.
\end{prp}

\begin{figure}[h]
\unitlength 0.1in
\begin{picture}(53.48,28.84)(13.40,-32.74)
%
\special{pn 8}%
\special{pa 1877 390}%
\special{pa 1370 1052}%
\special{fp}%
%
\special{pn 8}%
\special{pa 1374 1046}%
\special{pa 2359 1046}%
\special{fp}%
%
\special{pn 8}%
\special{pa 1877 390}%
\special{pa 2359 1052}%
\special{fp}%
%
\special{pn 8}%
\special{pa 1877 390}%
\special{pa 1877 833}%
\special{dt 0.045}%
\special{pa 1877 833}%
\special{pa 1877 832}%
\special{dt 0.045}%
\special{pa 1374 1046}%
\special{pa 1877 828}%
\special{dt 0.045}%
\special{pa 1877 828}%
\special{pa 1876 828}%
\special{dt 0.045}%
\special{pa 1877 828}%
\special{pa 2349 1046}%
\special{dt 0.045}%
\special{pa 2349 1046}%
\special{pa 2348 1046}%
\special{dt 0.045}%
%
\special{pn 20}%
\special{pa 1867 1348}%
\special{pa 1867 1156}%
\special{fp}%
\special{sh 1}%
\special{pa 1867 1156}%
\special{pa 1847 1223}%
\special{pa 1867 1209}%
\special{pa 1887 1223}%
\special{pa 1867 1156}%
\special{fp}%
%
\special{pn 8}%
\special{pa 1867 1386}%
\special{pa 1360 2048}%
\special{fp}%
%
\special{pn 8}%
\special{pa 1364 2043}%
\special{pa 2349 2043}%
\special{fp}%
%
\special{pn 8}%
\special{pa 1867 1386}%
\special{pa 2349 2048}%
\special{fp}%
%
\special{pn 8}%
\special{pa 1867 1386}%
\special{pa 1867 1829}%
\special{dt 0.045}%
\special{pa 1867 1829}%
\special{pa 1867 1828}%
\special{dt 0.045}%
\special{pa 1364 2043}%
\special{pa 1867 1824}%
\special{dt 0.045}%
\special{pa 1867 1824}%
\special{pa 1866 1824}%
\special{dt 0.045}%
\special{pa 1867 1824}%
\special{pa 2339 2043}%
\special{dt 0.045}%
\special{pa 2339 2043}%
\special{pa 2338 2043}%
\special{dt 0.045}%
%
\special{pn 8}%
\special{pa 1763 1517}%
\special{pa 1960 1517}%
\special{fp}%
%
\special{pn 8}%
\special{pa 1542 1807}%
\special{pa 2172 1807}%
\special{fp}%
%
\special{pn 8}%
\special{pa 1769 1517}%
\special{pa 1769 1867}%
\special{dt 0.045}%
\special{pa 1769 1867}%
\special{pa 1769 1866}%
\special{dt 0.045}%
\special{pa 1769 1867}%
\special{pa 1960 1867}%
\special{dt 0.045}%
\special{pa 1960 1867}%
\special{pa 1959 1867}%
\special{dt 0.045}%
\special{pa 1960 1522}%
\special{pa 1960 1867}%
\special{dt 0.045}%
\special{pa 1960 1867}%
\special{pa 1960 1866}%
\special{dt 0.045}%
%
\special{pn 8}%
\special{pa 1547 1807}%
\special{pa 1562 1950}%
\special{dt 0.045}%
\special{pa 1562 1950}%
\special{pa 1562 1949}%
\special{dt 0.045}%
\special{pa 1562 1955}%
\special{pa 2148 1955}%
\special{dt 0.045}%
\special{pa 2148 1955}%
\special{pa 2147 1955}%
\special{dt 0.045}%
\special{pa 2148 1955}%
\special{pa 2167 1802}%
\special{dt 0.045}%
\special{pa 2167 1802}%
\special{pa 2167 1803}%
\special{dt 0.045}%
%
\special{pn 20}%
\special{pa 1864 3274}%
\special{pa 1864 3082}%
\special{fp}%
\special{sh 1}%
\special{pa 1864 3082}%
\special{pa 1844 3149}%
\special{pa 1864 3135}%
\special{pa 1884 3149}%
\special{pa 1864 3082}%
\special{fp}%
%
\special{pn 8}%
\special{pa 1872 2371}%
\special{pa 1364 3033}%
\special{fp}%
%
\special{pn 8}%
\special{pa 1370 3027}%
\special{pa 2354 3027}%
\special{fp}%
%
\special{pn 8}%
\special{pa 1872 2371}%
\special{pa 2354 3033}%
\special{fp}%
%
\special{pn 8}%
\special{pa 1872 2371}%
\special{pa 1872 2814}%
\special{dt 0.045}%
\special{pa 1872 2814}%
\special{pa 1872 2813}%
\special{dt 0.045}%
\special{pa 1370 3027}%
\special{pa 1872 2809}%
\special{dt 0.045}%
\special{pa 1872 2809}%
\special{pa 1871 2809}%
\special{dt 0.045}%
\special{pa 1872 2809}%
\special{pa 2345 3027}%
\special{dt 0.045}%
\special{pa 2345 3027}%
\special{pa 2344 3027}%
\special{dt 0.045}%
%
\special{pn 8}%
\special{pa 1661 2639}%
\special{pa 2064 2639}%
\special{fp}%
%
\special{pn 8}%
\special{pa 1666 2639}%
\special{pa 1666 2896}%
\special{dt 0.045}%
\special{pa 1666 2896}%
\special{pa 1666 2895}%
\special{dt 0.045}%
\special{pa 1666 2896}%
\special{pa 2062 2896}%
\special{dt 0.045}%
\special{pa 2062 2896}%
\special{pa 2061 2896}%
\special{dt 0.045}%
\special{pa 2062 2896}%
\special{pa 2062 2644}%
\special{dt 0.045}%
\special{pa 2062 2644}%
\special{pa 2062 2645}%
\special{dt 0.045}%
%
\special{pn 8}%
\special{pa 3047 390}%
\special{pa 2539 1052}%
\special{fp}%
%
\special{pn 8}%
\special{pa 2544 1046}%
\special{pa 3529 1046}%
\special{fp}%
%
\special{pn 8}%
\special{pa 3047 390}%
\special{pa 3529 1052}%
\special{fp}%
%
\special{pn 8}%
\special{pa 3047 390}%
\special{pa 3047 833}%
\special{dt 0.045}%
\special{pa 3047 833}%
\special{pa 3047 832}%
\special{dt 0.045}%
\special{pa 2544 1046}%
\special{pa 3047 828}%
\special{dt 0.045}%
\special{pa 3047 828}%
\special{pa 3046 828}%
\special{dt 0.045}%
\special{pa 3047 828}%
\special{pa 3519 1046}%
\special{dt 0.045}%
\special{pa 3519 1046}%
\special{pa 3518 1046}%
\special{dt 0.045}%
%
\special{pn 8}%
\special{pa 2835 659}%
\special{pa 3238 659}%
\special{fp}%
%
\special{pn 8}%
\special{pa 2840 659}%
\special{pa 2840 915}%
\special{dt 0.045}%
\special{pa 2840 915}%
\special{pa 2840 914}%
\special{dt 0.045}%
\special{pa 2840 915}%
\special{pa 3236 915}%
\special{dt 0.045}%
\special{pa 3236 915}%
\special{pa 3235 915}%
\special{dt 0.045}%
\special{pa 3236 915}%
\special{pa 3236 664}%
\special{dt 0.045}%
\special{pa 3236 664}%
\special{pa 3236 665}%
\special{dt 0.045}%
%
\special{pn 20}%
\special{pa 1864 2103}%
\special{pa 1864 2289}%
\special{fp}%
\special{sh 1}%
\special{pa 1864 2289}%
\special{pa 1884 2222}%
\special{pa 1864 2236}%
\special{pa 1844 2222}%
\special{pa 1864 2289}%
\special{fp}%
%
\special{pn 8}%
\special{pa 3043 483}%
\special{pa 2845 741}%
\special{dt 0.045}%
\special{pa 2845 741}%
\special{pa 2845 740}%
\special{dt 0.045}%
\special{pa 2845 741}%
\special{pa 3232 741}%
\special{dt 0.045}%
\special{pa 3232 741}%
\special{pa 3231 741}%
\special{dt 0.045}%
\special{pa 3232 741}%
\special{pa 3043 483}%
\special{dt 0.045}%
\special{pa 3043 483}%
\special{pa 3043 484}%
\special{dt 0.045}%
%
\special{pn 8}%
\special{pa 3043 779}%
\special{pa 2845 866}%
\special{dt 0.045}%
\special{pa 2845 866}%
\special{pa 2846 866}%
\special{dt 0.045}%
\special{pa 2845 866}%
\special{pa 3232 866}%
\special{dt 0.045}%
\special{pa 3232 866}%
\special{pa 3231 866}%
\special{dt 0.045}%
\special{pa 3232 866}%
\special{pa 3043 784}%
\special{dt 0.045}%
\special{pa 3043 784}%
\special{pa 3044 784}%
\special{dt 0.045}%
%
\special{pn 20}%
\special{pa 3016 1140}%
\special{pa 3016 1326}%
\special{fp}%
\special{sh 1}%
\special{pa 3016 1326}%
\special{pa 3036 1259}%
\special{pa 3016 1273}%
\special{pa 2996 1259}%
\special{pa 3016 1326}%
\special{fp}%
%
\special{pn 8}%
\special{pa 3038 1391}%
\special{pa 2530 2053}%
\special{fp}%
%
\special{pn 8}%
\special{pa 2535 2048}%
\special{pa 3520 2048}%
\special{fp}%
%
\special{pn 8}%
\special{pa 3038 1391}%
\special{pa 3520 2053}%
\special{fp}%
%
\special{pn 8}%
\special{pa 3038 1391}%
\special{pa 3038 1835}%
\special{dt 0.045}%
\special{pa 3038 1835}%
\special{pa 3038 1834}%
\special{dt 0.045}%
\special{pa 2535 2048}%
\special{pa 3038 1829}%
\special{dt 0.045}%
\special{pa 3038 1829}%
\special{pa 3037 1829}%
\special{dt 0.045}%
\special{pa 3038 1829}%
\special{pa 3510 2048}%
\special{dt 0.045}%
\special{pa 3510 2048}%
\special{pa 3509 2048}%
\special{dt 0.045}%
%
\special{pn 8}%
\special{pa 2826 1659}%
\special{pa 3229 1659}%
\special{fp}%
%
\special{pn 8}%
\special{pa 2832 1659}%
\special{pa 2832 1917}%
\special{dt 0.045}%
\special{pa 2832 1917}%
\special{pa 2832 1916}%
\special{dt 0.045}%
\special{pa 2832 1917}%
\special{pa 3228 1917}%
\special{dt 0.045}%
\special{pa 3228 1917}%
\special{pa 3227 1917}%
\special{dt 0.045}%
\special{pa 3228 1917}%
\special{pa 3228 1665}%
\special{dt 0.045}%
\special{pa 3228 1665}%
\special{pa 3228 1666}%
\special{dt 0.045}%
%
\special{pn 8}%
\special{pa 3034 1709}%
\special{pa 2836 1780}%
\special{dt 0.045}%
\special{pa 2836 1780}%
\special{pa 2837 1780}%
\special{dt 0.045}%
\special{pa 2836 1780}%
\special{pa 3223 1780}%
\special{dt 0.045}%
\special{pa 3223 1780}%
\special{pa 3222 1780}%
\special{dt 0.045}%
\special{pa 3223 1780}%
\special{pa 3034 1709}%
\special{dt 0.045}%
\special{pa 3034 1709}%
\special{pa 3035 1709}%
\special{dt 0.045}%
%
\special{pn 20}%
\special{pa 3025 2283}%
\special{pa 3025 2091}%
\special{fp}%
\special{sh 1}%
\special{pa 3025 2091}%
\special{pa 3005 2158}%
\special{pa 3025 2144}%
\special{pa 3045 2158}%
\special{pa 3025 2091}%
\special{fp}%
%
\special{pn 8}%
\special{pa 3038 2371}%
\special{pa 2530 3033}%
\special{fp}%
%
\special{pn 8}%
\special{pa 2535 3027}%
\special{pa 3520 3027}%
\special{fp}%
%
\special{pn 8}%
\special{pa 3038 2371}%
\special{pa 3520 3033}%
\special{fp}%
%
\special{pn 8}%
\special{pa 3038 2371}%
\special{pa 3038 2814}%
\special{dt 0.045}%
\special{pa 3038 2814}%
\special{pa 3038 2813}%
\special{dt 0.045}%
\special{pa 2535 3027}%
\special{pa 3038 2809}%
\special{dt 0.045}%
\special{pa 3038 2809}%
\special{pa 3037 2809}%
\special{dt 0.045}%
\special{pa 3038 2809}%
\special{pa 3510 3027}%
\special{dt 0.045}%
\special{pa 3510 3027}%
\special{pa 3509 3027}%
\special{dt 0.045}%
%
\special{pn 8}%
\special{pa 2826 2639}%
\special{pa 3229 2639}%
\special{fp}%
%
\special{pn 8}%
\special{pa 2832 2639}%
\special{pa 2832 2896}%
\special{dt 0.045}%
\special{pa 2832 2896}%
\special{pa 2832 2895}%
\special{dt 0.045}%
\special{pa 2832 2896}%
\special{pa 3228 2896}%
\special{dt 0.045}%
\special{pa 3228 2896}%
\special{pa 3227 2896}%
\special{dt 0.045}%
\special{pa 3228 2896}%
\special{pa 3228 2644}%
\special{dt 0.045}%
\special{pa 3228 2644}%
\special{pa 3228 2645}%
\special{dt 0.045}%
%
\special{pn 8}%
\special{pa 3034 2688}%
\special{pa 2836 2759}%
\special{dt 0.045}%
\special{pa 2836 2759}%
\special{pa 2837 2759}%
\special{dt 0.045}%
\special{pa 2836 2759}%
\special{pa 3223 2759}%
\special{dt 0.045}%
\special{pa 3223 2759}%
\special{pa 3222 2759}%
\special{dt 0.045}%
\special{pa 3223 2759}%
\special{pa 3034 2688}%
\special{dt 0.045}%
\special{pa 3034 2688}%
\special{pa 3035 2688}%
\special{dt 0.045}%
%
\special{pn 20}%
\special{pa 3025 3077}%
\special{pa 3025 3263}%
\special{fp}%
\special{sh 1}%
\special{pa 3025 3263}%
\special{pa 3045 3196}%
\special{pa 3025 3210}%
\special{pa 3005 3196}%
\special{pa 3025 3263}%
\special{fp}%
%
\special{pn 8}%
\special{pa 2953 2470}%
\special{pa 3070 2639}%
\special{dt 0.045}%
\special{pa 3070 2639}%
\special{pa 3070 2638}%
\special{dt 0.045}%
\special{pa 3070 2639}%
\special{pa 3070 2759}%
\special{dt 0.045}%
\special{pa 3070 2759}%
\special{pa 3070 2758}%
\special{dt 0.045}%
\special{pa 2962 2470}%
\special{pa 2962 2704}%
\special{dt 0.045}%
\special{pa 2962 2704}%
\special{pa 2962 2703}%
\special{dt 0.045}%
\special{pa 2962 2704}%
\special{pa 3070 2754}%
\special{dt 0.045}%
\special{pa 3070 2754}%
\special{pa 3069 2754}%
\special{dt 0.045}%
%
\special{pn 8}%
\special{pa 4325 395}%
\special{pa 3817 1058}%
\special{fp}%
%
\special{pn 8}%
\special{pa 3822 1052}%
\special{pa 4807 1052}%
\special{fp}%
%
\special{pn 8}%
\special{pa 4325 395}%
\special{pa 4807 1058}%
\special{fp}%
%
\special{pn 8}%
\special{pa 4325 395}%
\special{pa 4325 838}%
\special{dt 0.045}%
\special{pa 4325 838}%
\special{pa 4325 837}%
\special{dt 0.045}%
\special{pa 3822 1052}%
\special{pa 4325 833}%
\special{dt 0.045}%
\special{pa 4325 833}%
\special{pa 4324 833}%
\special{dt 0.045}%
\special{pa 4325 833}%
\special{pa 4797 1052}%
\special{dt 0.045}%
\special{pa 4797 1052}%
\special{pa 4796 1052}%
\special{dt 0.045}%
%
\special{pn 8}%
\special{pa 4113 664}%
\special{pa 4516 664}%
\special{fp}%
%
\special{pn 8}%
\special{pa 4118 664}%
\special{pa 4118 921}%
\special{dt 0.045}%
\special{pa 4118 921}%
\special{pa 4118 920}%
\special{dt 0.045}%
\special{pa 4118 921}%
\special{pa 4514 921}%
\special{dt 0.045}%
\special{pa 4514 921}%
\special{pa 4513 921}%
\special{dt 0.045}%
\special{pa 4514 921}%
\special{pa 4514 668}%
\special{dt 0.045}%
\special{pa 4514 668}%
\special{pa 4514 669}%
\special{dt 0.045}%
%
\special{pn 8}%
\special{pa 4321 713}%
\special{pa 4123 784}%
\special{dt 0.045}%
\special{pa 4123 784}%
\special{pa 4124 784}%
\special{dt 0.045}%
\special{pa 4123 784}%
\special{pa 4510 784}%
\special{dt 0.045}%
\special{pa 4510 784}%
\special{pa 4509 784}%
\special{dt 0.045}%
\special{pa 4510 784}%
\special{pa 4321 713}%
\special{dt 0.045}%
\special{pa 4321 713}%
\special{pa 4322 713}%
\special{dt 0.045}%
%
\special{pn 20}%
\special{pa 4321 1320}%
\special{pa 4321 1129}%
\special{fp}%
\special{sh 1}%
\special{pa 4321 1129}%
\special{pa 4301 1196}%
\special{pa 4321 1182}%
\special{pa 4341 1196}%
\special{pa 4321 1129}%
\special{fp}%
%
\special{pn 20}%
\special{pa 4321 2097}%
\special{pa 4321 2283}%
\special{fp}%
\special{sh 1}%
\special{pa 4321 2283}%
\special{pa 4341 2216}%
\special{pa 4321 2230}%
\special{pa 4301 2216}%
\special{pa 4321 2283}%
\special{fp}%
%
\special{pn 8}%
\special{pa 4334 2365}%
\special{pa 3826 3027}%
\special{fp}%
%
\special{pn 8}%
\special{pa 3831 3022}%
\special{pa 4816 3022}%
\special{fp}%
%
\special{pn 8}%
\special{pa 4334 2365}%
\special{pa 4816 3027}%
\special{fp}%
%
\special{pn 8}%
\special{pa 4334 2365}%
\special{pa 4334 2809}%
\special{dt 0.045}%
\special{pa 4334 2809}%
\special{pa 4334 2808}%
\special{dt 0.045}%
\special{pa 3831 3022}%
\special{pa 4334 2803}%
\special{dt 0.045}%
\special{pa 4334 2803}%
\special{pa 4333 2803}%
\special{dt 0.045}%
\special{pa 4334 2803}%
\special{pa 4806 3022}%
\special{dt 0.045}%
\special{pa 4806 3022}%
\special{pa 4805 3022}%
\special{dt 0.045}%
%
\special{pn 8}%
\special{pa 4122 2634}%
\special{pa 4525 2634}%
\special{fp}%
%
\special{pn 8}%
\special{pa 4128 2634}%
\special{pa 4128 2890}%
\special{dt 0.045}%
\special{pa 4128 2890}%
\special{pa 4128 2889}%
\special{dt 0.045}%
\special{pa 4128 2890}%
\special{pa 4524 2890}%
\special{dt 0.045}%
\special{pa 4524 2890}%
\special{pa 4523 2890}%
\special{dt 0.045}%
\special{pa 4524 2890}%
\special{pa 4524 2639}%
\special{dt 0.045}%
\special{pa 4524 2639}%
\special{pa 4524 2640}%
\special{dt 0.045}%
%
\special{pn 8}%
\special{pa 4330 2683}%
\special{pa 4132 2754}%
\special{dt 0.045}%
\special{pa 4132 2754}%
\special{pa 4133 2754}%
\special{dt 0.045}%
\special{pa 4132 2754}%
\special{pa 4519 2754}%
\special{dt 0.045}%
\special{pa 4519 2754}%
\special{pa 4518 2754}%
\special{dt 0.045}%
\special{pa 4519 2754}%
\special{pa 4330 2683}%
\special{dt 0.045}%
\special{pa 4330 2683}%
\special{pa 4331 2683}%
\special{dt 0.045}%
%
\special{pn 20}%
\special{pa 4330 3263}%
\special{pa 4330 3071}%
\special{fp}%
\special{sh 1}%
\special{pa 4330 3071}%
\special{pa 4310 3138}%
\special{pa 4330 3124}%
\special{pa 4350 3138}%
\special{pa 4330 3071}%
\special{fp}%
%
\special{pn 20}%
\special{pa 5626 1129}%
\special{pa 5626 1315}%
\special{fp}%
\special{sh 1}%
\special{pa 5626 1315}%
\special{pa 5646 1248}%
\special{pa 5626 1262}%
\special{pa 5606 1248}%
\special{pa 5626 1315}%
\special{fp}%
%
\special{pn 8}%
\special{pa 5630 1391}%
\special{pa 5122 2053}%
\special{fp}%
%
\special{pn 8}%
\special{pa 5127 2048}%
\special{pa 6112 2048}%
\special{fp}%
%
\special{pn 8}%
\special{pa 5630 1391}%
\special{pa 6112 2053}%
\special{fp}%
%
\special{pn 8}%
\special{pa 5630 1391}%
\special{pa 5630 1835}%
\special{dt 0.045}%
\special{pa 5630 1835}%
\special{pa 5630 1834}%
\special{dt 0.045}%
\special{pa 5127 2048}%
\special{pa 5630 1829}%
\special{dt 0.045}%
\special{pa 5630 1829}%
\special{pa 5629 1829}%
\special{dt 0.045}%
\special{pa 5630 1829}%
\special{pa 6102 2048}%
\special{dt 0.045}%
\special{pa 6102 2048}%
\special{pa 6101 2048}%
\special{dt 0.045}%
%
\special{pn 8}%
\special{pa 5418 1659}%
\special{pa 5821 1659}%
\special{fp}%
%
\special{pn 8}%
\special{pa 5424 1659}%
\special{pa 5424 1917}%
\special{dt 0.045}%
\special{pa 5424 1917}%
\special{pa 5424 1916}%
\special{dt 0.045}%
\special{pa 5424 1917}%
\special{pa 5820 1917}%
\special{dt 0.045}%
\special{pa 5820 1917}%
\special{pa 5819 1917}%
\special{dt 0.045}%
\special{pa 5820 1917}%
\special{pa 5820 1665}%
\special{dt 0.045}%
\special{pa 5820 1665}%
\special{pa 5820 1666}%
\special{dt 0.045}%
%
\special{pn 8}%
\special{pa 5626 1709}%
\special{pa 5428 1780}%
\special{dt 0.045}%
\special{pa 5428 1780}%
\special{pa 5429 1780}%
\special{dt 0.045}%
\special{pa 5428 1780}%
\special{pa 5815 1780}%
\special{dt 0.045}%
\special{pa 5815 1780}%
\special{pa 5814 1780}%
\special{dt 0.045}%
\special{pa 5815 1780}%
\special{pa 5626 1709}%
\special{dt 0.045}%
\special{pa 5626 1709}%
\special{pa 5627 1709}%
\special{dt 0.045}%
%
\special{pn 20}%
\special{pa 5626 3055}%
\special{pa 5626 3241}%
\special{fp}%
\special{sh 1}%
\special{pa 5626 3241}%
\special{pa 5646 3174}%
\special{pa 5626 3188}%
\special{pa 5606 3174}%
\special{pa 5626 3241}%
\special{fp}%
\put(55.7000,-34.2000){\makebox(0,0)[lb]{$\mathcal S$}}%
%
\special{pn 8}%
\special{pa 3080 2641}%
\special{pa 2960 2483}%
\special{fp}%
%
\special{pn 8}%
\special{pa 5638 397}%
\special{pa 5130 1059}%
\special{fp}%
%
\special{pn 8}%
\special{pa 5135 1053}%
\special{pa 6120 1053}%
\special{fp}%
%
\special{pn 8}%
\special{pa 5638 397}%
\special{pa 6120 1059}%
\special{fp}%
%
\special{pn 8}%
\special{pa 5638 397}%
\special{pa 5638 839}%
\special{dt 0.045}%
\special{pa 5638 839}%
\special{pa 5638 838}%
\special{dt 0.045}%
\special{pa 5135 1053}%
\special{pa 5638 835}%
\special{dt 0.045}%
\special{pa 5638 835}%
\special{pa 5637 835}%
\special{dt 0.045}%
\special{pa 5638 835}%
\special{pa 6110 1053}%
\special{dt 0.045}%
\special{pa 6110 1053}%
\special{pa 6109 1053}%
\special{dt 0.045}%
%
\special{pn 8}%
\special{pa 5426 666}%
\special{pa 5829 666}%
\special{fp}%
%
\special{pn 8}%
\special{pa 5432 666}%
\special{pa 5432 922}%
\special{dt 0.045}%
\special{pa 5432 922}%
\special{pa 5432 921}%
\special{dt 0.045}%
\special{pa 5432 922}%
\special{pa 5828 922}%
\special{dt 0.045}%
\special{pa 5828 922}%
\special{pa 5827 922}%
\special{dt 0.045}%
\special{pa 5828 922}%
\special{pa 5828 670}%
\special{dt 0.045}%
\special{pa 5828 670}%
\special{pa 5828 671}%
\special{dt 0.045}%
%
\special{pn 8}%
\special{pa 5634 714}%
\special{pa 5436 785}%
\special{dt 0.045}%
\special{pa 5436 785}%
\special{pa 5437 785}%
\special{dt 0.045}%
\special{pa 5436 785}%
\special{pa 5823 785}%
\special{dt 0.045}%
\special{pa 5823 785}%
\special{pa 5822 785}%
\special{dt 0.045}%
\special{pa 5823 785}%
\special{pa 5634 714}%
\special{dt 0.045}%
\special{pa 5634 714}%
\special{pa 5635 714}%
\special{dt 0.045}%
%
\special{pn 8}%
\special{pa 5688 736}%
\special{pa 5544 785}%
\special{dt 0.045}%
\special{pa 5544 785}%
\special{pa 5545 785}%
\special{dt 0.045}%
\special{pa 6678 873}%
\special{pa 6678 873}%
\special{dt 0.045}%
%
\special{pn 8}%
\special{pa 5553 785}%
\special{pa 5553 922}%
\special{dt 0.045}%
\special{pa 5553 922}%
\special{pa 5553 921}%
\special{dt 0.045}%
\special{pa 5553 922}%
\special{pa 5697 862}%
\special{dt 0.045}%
\special{pa 5697 862}%
\special{pa 5696 862}%
\special{dt 0.045}%
\special{pa 5697 862}%
\special{pa 5697 742}%
\special{dt 0.045}%
\special{pa 5697 742}%
\special{pa 5697 743}%
\special{dt 0.045}%
\put(19.0000,-13.8000){\makebox(0,0)[lb]{Step 1}}%
\put(19.0000,-23.0000){\makebox(0,0)[lb]{Step 2}}%
\put(19.0000,-33.0000){\makebox(0,0)[lb]{Step 3}}%
\put(30.7000,-13.5000){\makebox(0,0)[lb]{Step 4}}%
\put(30.8000,-22.8000){\makebox(0,0)[lb]{Step 5}}%
\put(30.9000,-32.7000){\makebox(0,0)[lb]{Step 6}}%
\put(43.7000,-13.3000){\makebox(0,0)[lb]{Step 7}}%
\put(43.8000,-22.7000){\makebox(0,0)[lb]{Step 8}}%
\put(43.8000,-32.5000){\makebox(0,0)[lb]{Step 9}}%
\put(56.8000,-13.3000){\makebox(0,0)[lb]{Step 10}}%
\put(56.9000,-22.7000){\makebox(0,0)[lb]{Step 11}}%
\put(56.9000,-32.3000){\makebox(0,0)[lb]{Step 12}}%
%
\special{pn 20}%
\special{pa 1890 410}%
\special{pa 1890 820}%
\special{fp}%
%
\special{pn 20}%
\special{pa 1380 1040}%
\special{pa 2350 1040}%
\special{fp}%
%
\special{pn 8}%
\special{pa 1770 1730}%
\special{pa 1690 1810}%
\special{fp}%
\special{pa 1770 1670}%
\special{pa 1630 1810}%
\special{fp}%
\special{pa 1770 1610}%
\special{pa 1570 1810}%
\special{fp}%
\special{pa 1770 1550}%
\special{pa 1680 1640}%
\special{fp}%
%
\special{pn 4}%
\special{pa 1870 1690}%
\special{pa 1770 1790}%
\special{fp}%
\special{pa 1870 1750}%
\special{pa 1810 1810}%
\special{fp}%
\special{pa 1870 1630}%
\special{pa 1770 1730}%
\special{fp}%
\special{pa 1870 1570}%
\special{pa 1770 1670}%
\special{fp}%
\special{pa 1860 1520}%
\special{pa 1770 1610}%
\special{fp}%
\special{pa 1800 1520}%
\special{pa 1770 1550}%
\special{fp}%
%
\special{pn 4}%
\special{pa 1960 1660}%
\special{pa 1870 1750}%
\special{fp}%
\special{pa 1960 1720}%
\special{pa 1880 1800}%
\special{fp}%
\special{pa 1960 1780}%
\special{pa 1930 1810}%
\special{fp}%
\special{pa 1960 1600}%
\special{pa 1870 1690}%
\special{fp}%
\special{pa 1960 1540}%
\special{pa 1870 1630}%
\special{fp}%
\special{pa 1920 1520}%
\special{pa 1870 1570}%
\special{fp}%
%
\special{pn 4}%
\special{pa 2070 1670}%
\special{pa 1960 1780}%
\special{fp}%
\special{pa 2100 1700}%
\special{pa 1990 1810}%
\special{fp}%
\special{pa 2120 1740}%
\special{pa 2050 1810}%
\special{fp}%
\special{pa 2150 1770}%
\special{pa 2110 1810}%
\special{fp}%
\special{pa 2040 1640}%
\special{pa 1960 1720}%
\special{fp}%
\special{pa 2020 1600}%
\special{pa 1960 1660}%
\special{fp}%
\special{pa 1990 1570}%
\special{pa 1960 1600}%
\special{fp}%
%
\special{pn 4}%
\special{pa 1630 1810}%
\special{pa 1560 1880}%
\special{fp}%
\special{pa 1690 1810}%
\special{pa 1560 1940}%
\special{fp}%
\special{pa 1750 1810}%
\special{pa 1650 1910}%
\special{fp}%
%
\special{pn 4}%
\special{pa 1930 1870}%
\special{pa 1840 1960}%
\special{fp}%
\special{pa 1870 1870}%
\special{pa 1780 1960}%
\special{fp}%
\special{pa 1810 1870}%
\special{pa 1720 1960}%
\special{fp}%
\special{pa 1740 1880}%
\special{pa 1660 1960}%
\special{fp}%
\special{pa 1980 1880}%
\special{pa 1900 1960}%
\special{fp}%
\special{pa 2030 1890}%
\special{pa 1960 1960}%
\special{fp}%
\special{pa 2070 1910}%
\special{pa 2020 1960}%
\special{fp}%
%
\special{pn 4}%
\special{pa 1900 1840}%
\special{pa 1870 1870}%
\special{fp}%
%
\special{pn 4}%
\special{pa 2160 1820}%
\special{pa 2070 1910}%
\special{fp}%
\special{pa 2160 1880}%
\special{pa 2110 1930}%
\special{fp}%
\special{pa 2110 1810}%
\special{pa 2030 1890}%
\special{fp}%
\special{pa 2050 1810}%
\special{pa 1990 1870}%
\special{fp}%
\special{pa 1990 1810}%
\special{pa 1960 1840}%
\special{fp}%
%
\special{pn 4}%
\special{pa 1870 2470}%
\special{pa 1700 2640}%
\special{fp}%
\special{pa 1870 2530}%
\special{pa 1760 2640}%
\special{fp}%
\special{pa 1870 2590}%
\special{pa 1820 2640}%
\special{fp}%
\special{pa 1870 2410}%
\special{pa 1750 2530}%
\special{fp}%
%
\special{pn 4}%
\special{pa 1990 2530}%
\special{pa 1880 2640}%
\special{fp}%
\special{pa 1960 2500}%
\special{pa 1870 2590}%
\special{fp}%
\special{pa 1940 2460}%
\special{pa 1870 2530}%
\special{fp}%
\special{pa 1910 2430}%
\special{pa 1870 2470}%
\special{fp}%
\special{pa 2010 2570}%
\special{pa 1940 2640}%
\special{fp}%
\special{pa 2040 2600}%
\special{pa 2000 2640}%
\special{fp}%
%
\special{pn 4}%
\special{pa 1880 2820}%
\special{pa 1800 2900}%
\special{fp}%
\special{pa 1790 2850}%
\special{pa 1740 2900}%
\special{fp}%
\special{pa 1920 2840}%
\special{pa 1860 2900}%
\special{fp}%
\special{pa 1960 2860}%
\special{pa 1920 2900}%
\special{fp}%
%
\special{pn 4}%
\special{pa 3000 740}%
\special{pa 2890 850}%
\special{fp}%
\special{pa 2940 740}%
\special{pa 2840 840}%
\special{fp}%
\special{pa 2880 740}%
\special{pa 2840 780}%
\special{fp}%
\special{pa 3050 750}%
\special{pa 3010 790}%
\special{fp}%
%
\special{pn 4}%
\special{pa 3230 750}%
\special{pa 3150 830}%
\special{fp}%
\special{pa 3240 800}%
\special{pa 3190 850}%
\special{fp}%
\special{pa 3180 740}%
\special{pa 3110 810}%
\special{fp}%
\special{pa 3120 740}%
\special{pa 3070 790}%
\special{fp}%
%
\special{pn 4}%
\special{pa 3020 660}%
\special{pa 2940 740}%
\special{fp}%
\special{pa 2960 660}%
\special{pa 2880 740}%
\special{fp}%
\special{pa 3050 690}%
\special{pa 3000 740}%
\special{fp}%
%
\special{pn 4}%
\special{pa 3140 660}%
\special{pa 3060 740}%
\special{fp}%
\special{pa 3180 680}%
\special{pa 3120 740}%
\special{fp}%
\special{pa 3210 710}%
\special{pa 3180 740}%
\special{fp}%
\special{pa 3080 660}%
\special{pa 3050 690}%
\special{fp}%
%
\special{pn 4}%
\special{pa 3050 570}%
\special{pa 2960 660}%
\special{fp}%
\special{pa 3050 630}%
\special{pa 3020 660}%
\special{fp}%
\special{pa 3050 510}%
\special{pa 2940 620}%
\special{fp}%
%
\special{pn 4}%
\special{pa 3130 610}%
\special{pa 3080 660}%
\special{fp}%
\special{pa 3110 570}%
\special{pa 3050 630}%
\special{fp}%
\special{pa 3080 540}%
\special{pa 3050 570}%
\special{fp}%
%
\special{pn 4}%
\special{pa 3110 810}%
\special{pa 3080 840}%
\special{fp}%
\special{pa 3150 830}%
\special{pa 3120 860}%
\special{fp}%
%
\special{pn 4}%
\special{pa 3000 800}%
\special{pa 2930 870}%
\special{fp}%
%
\special{pn 4}%
\special{pa 3150 1550}%
\special{pa 3050 1650}%
\special{fp}%
\special{pa 3180 1580}%
\special{pa 3100 1660}%
\special{fp}%
\special{pa 3200 1620}%
\special{pa 3160 1660}%
\special{fp}%
\special{pa 3130 1510}%
\special{pa 3040 1600}%
\special{fp}%
\special{pa 3100 1480}%
\special{pa 3040 1540}%
\special{fp}%
\special{pa 3080 1440}%
\special{pa 3040 1480}%
\special{fp}%
%
\special{pn 4}%
\special{pa 3220 1660}%
\special{pa 3130 1750}%
\special{fp}%
\special{pa 3230 1710}%
\special{pa 3180 1760}%
\special{fp}%
\special{pa 3160 1660}%
\special{pa 3090 1730}%
\special{fp}%
\special{pa 3100 1660}%
\special{pa 3050 1710}%
\special{fp}%
%
\special{pn 4}%
\special{pa 2960 2580}%
\special{pa 2900 2640}%
\special{fp}%
\special{pa 2960 2520}%
\special{pa 2850 2630}%
\special{fp}%
%
\special{pn 4}%
\special{pa 3020 2580}%
\special{pa 2970 2630}%
\special{fp}%
\special{pa 3000 2540}%
\special{pa 2960 2580}%
\special{fp}%
%
\special{pn 4}%
\special{pa 2950 2650}%
\special{pa 2840 2760}%
\special{fp}%
\special{pa 2900 2640}%
\special{pa 2830 2710}%
\special{fp}%
%
\special{pn 4}%
\special{pa 3020 2640}%
\special{pa 2970 2690}%
\special{fp}%
%
\special{pn 4}%
\special{pa 3000 2720}%
\special{pa 2960 2760}%
\special{fp}%
%
\special{pn 8}%
\special{pa 4338 1370}%
\special{pa 3830 2032}%
\special{fp}%
%
\special{pn 8}%
\special{pa 3835 2026}%
\special{pa 4820 2026}%
\special{fp}%
%
\special{pn 8}%
\special{pa 4338 1370}%
\special{pa 4820 2032}%
\special{fp}%
%
\special{pn 8}%
\special{pa 4338 1370}%
\special{pa 4338 1813}%
\special{dt 0.045}%
\special{pa 4338 1813}%
\special{pa 4338 1812}%
\special{dt 0.045}%
\special{pa 3835 2026}%
\special{pa 4338 1808}%
\special{dt 0.045}%
\special{pa 4338 1808}%
\special{pa 4337 1808}%
\special{dt 0.045}%
\special{pa 4338 1808}%
\special{pa 4810 2026}%
\special{dt 0.045}%
\special{pa 4810 2026}%
\special{pa 4809 2026}%
\special{dt 0.045}%
%
\special{pn 8}%
\special{pa 4126 1638}%
\special{pa 4529 1638}%
\special{fp}%
%
\special{pn 8}%
\special{pa 4132 1638}%
\special{pa 4132 1895}%
\special{dt 0.045}%
\special{pa 4132 1895}%
\special{pa 4132 1894}%
\special{dt 0.045}%
\special{pa 4132 1895}%
\special{pa 4528 1895}%
\special{dt 0.045}%
\special{pa 4528 1895}%
\special{pa 4527 1895}%
\special{dt 0.045}%
\special{pa 4528 1895}%
\special{pa 4528 1643}%
\special{dt 0.045}%
\special{pa 4528 1643}%
\special{pa 4528 1644}%
\special{dt 0.045}%
%
\special{pn 8}%
\special{pa 4334 1687}%
\special{pa 4136 1758}%
\special{dt 0.045}%
\special{pa 4136 1758}%
\special{pa 4137 1758}%
\special{dt 0.045}%
\special{pa 4136 1758}%
\special{pa 4523 1758}%
\special{dt 0.045}%
\special{pa 4523 1758}%
\special{pa 4522 1758}%
\special{dt 0.045}%
\special{pa 4523 1758}%
\special{pa 4334 1687}%
\special{dt 0.045}%
\special{pa 4334 1687}%
\special{pa 4335 1687}%
\special{dt 0.045}%
%
\special{pn 8}%
\special{pa 4253 1469}%
\special{pa 4370 1638}%
\special{dt 0.045}%
\special{pa 4370 1638}%
\special{pa 4370 1637}%
\special{dt 0.045}%
\special{pa 4370 1638}%
\special{pa 4370 1758}%
\special{dt 0.045}%
\special{pa 4370 1758}%
\special{pa 4370 1757}%
\special{dt 0.045}%
\special{pa 4262 1469}%
\special{pa 4262 1703}%
\special{dt 0.045}%
\special{pa 4262 1703}%
\special{pa 4262 1702}%
\special{dt 0.045}%
\special{pa 4262 1703}%
\special{pa 4370 1753}%
\special{dt 0.045}%
\special{pa 4370 1753}%
\special{pa 4369 1753}%
\special{dt 0.045}%
%
\special{pn 8}%
\special{pa 4380 1640}%
\special{pa 4260 1482}%
\special{fp}%
%
\special{pn 4}%
\special{pa 4260 1579}%
\special{pa 4200 1639}%
\special{fp}%
\special{pa 4260 1519}%
\special{pa 4150 1629}%
\special{fp}%
%
\special{pn 4}%
\special{pa 4320 1579}%
\special{pa 4270 1629}%
\special{fp}%
\special{pa 4300 1539}%
\special{pa 4260 1579}%
\special{fp}%
%
\special{pn 4}%
\special{pa 4250 1649}%
\special{pa 4140 1759}%
\special{fp}%
\special{pa 4200 1639}%
\special{pa 4130 1709}%
\special{fp}%
%
\special{pn 4}%
\special{pa 4320 1639}%
\special{pa 4270 1689}%
\special{fp}%
%
\special{pn 4}%
\special{pa 4300 1719}%
\special{pa 4260 1759}%
\special{fp}%
%
\special{pn 4}%
\special{pa 4230 2750}%
\special{pa 4130 2850}%
\special{fp}%
\special{pa 4290 2750}%
\special{pa 4170 2870}%
\special{fp}%
\special{pa 4330 2770}%
\special{pa 4280 2820}%
\special{fp}%
\special{pa 4170 2750}%
\special{pa 4130 2790}%
\special{fp}%
%
\special{pn 4}%
\special{pa 4330 2710}%
\special{pa 4290 2750}%
\special{fp}%
\special{pa 4280 2700}%
\special{pa 4230 2750}%
\special{fp}%
%
\special{pn 4}%
\special{pa 5820 800}%
\special{pa 5740 880}%
\special{fp}%
\special{pa 5830 850}%
\special{pa 5780 900}%
\special{fp}%
\special{pa 5770 790}%
\special{pa 5700 860}%
\special{fp}%
%
\special{pn 4}%
\special{pa 5690 870}%
\special{pa 5640 920}%
\special{fp}%
\special{pa 5740 880}%
\special{pa 5700 920}%
\special{fp}%
%
\special{pn 4}%
\special{pa 5640 800}%
\special{pa 5570 870}%
\special{fp}%
\special{pa 5590 790}%
\special{pa 5550 830}%
\special{fp}%
%
\special{pn 4}%
\special{pa 5700 800}%
\special{pa 5660 840}%
\special{fp}%
%
\special{pn 4}%
\special{pa 5690 750}%
\special{pa 5650 790}%
\special{fp}%
%
\special{pn 4}%
\special{pa 5740 760}%
\special{pa 5710 790}%
\special{fp}%
%
\special{pn 4}%
\special{pa 5560 1780}%
\special{pa 5430 1910}%
\special{fp}%
\special{pa 5500 1780}%
\special{pa 5420 1860}%
\special{fp}%
\special{pa 5620 1780}%
\special{pa 5530 1870}%
\special{fp}%
%
\special{pn 4}%
\special{pa 5620 1720}%
\special{pa 5560 1780}%
\special{fp}%
%
\special{pn 8}%
\special{pa 5648 2350}%
\special{pa 5140 3012}%
\special{fp}%
%
\special{pn 8}%
\special{pa 5145 3006}%
\special{pa 6130 3006}%
\special{fp}%
%
\special{pn 8}%
\special{pa 5648 2350}%
\special{pa 6130 3012}%
\special{fp}%
%
\special{pn 8}%
\special{pa 5648 2350}%
\special{pa 5648 2792}%
\special{dt 0.045}%
\special{pa 5648 2792}%
\special{pa 5648 2791}%
\special{dt 0.045}%
\special{pa 5145 3006}%
\special{pa 5648 2788}%
\special{dt 0.045}%
\special{pa 5648 2788}%
\special{pa 5647 2788}%
\special{dt 0.045}%
\special{pa 5648 2788}%
\special{pa 6120 3006}%
\special{dt 0.045}%
\special{pa 6120 3006}%
\special{pa 6119 3006}%
\special{dt 0.045}%
%
\special{pn 8}%
\special{pa 5436 2619}%
\special{pa 5839 2619}%
\special{fp}%
%
\special{pn 8}%
\special{pa 5442 2619}%
\special{pa 5442 2875}%
\special{dt 0.045}%
\special{pa 5442 2875}%
\special{pa 5442 2874}%
\special{dt 0.045}%
\special{pa 5442 2875}%
\special{pa 5838 2875}%
\special{dt 0.045}%
\special{pa 5838 2875}%
\special{pa 5837 2875}%
\special{dt 0.045}%
\special{pa 5838 2875}%
\special{pa 5838 2623}%
\special{dt 0.045}%
\special{pa 5838 2623}%
\special{pa 5838 2624}%
\special{dt 0.045}%
%
\special{pn 8}%
\special{pa 5644 2667}%
\special{pa 5446 2738}%
\special{dt 0.045}%
\special{pa 5446 2738}%
\special{pa 5447 2738}%
\special{dt 0.045}%
\special{pa 5446 2738}%
\special{pa 5833 2738}%
\special{dt 0.045}%
\special{pa 5833 2738}%
\special{pa 5832 2738}%
\special{dt 0.045}%
\special{pa 5833 2738}%
\special{pa 5644 2667}%
\special{dt 0.045}%
\special{pa 5644 2667}%
\special{pa 5645 2667}%
\special{dt 0.045}%
%
\special{pn 8}%
\special{pa 5698 2689}%
\special{pa 5554 2738}%
\special{dt 0.045}%
\special{pa 5554 2738}%
\special{pa 5555 2738}%
\special{dt 0.045}%
\special{pa 6688 2826}%
\special{pa 6688 2826}%
\special{dt 0.045}%
%
\special{pn 8}%
\special{pa 5563 2738}%
\special{pa 5563 2875}%
\special{dt 0.045}%
\special{pa 5563 2875}%
\special{pa 5563 2874}%
\special{dt 0.045}%
\special{pa 5563 2875}%
\special{pa 5707 2815}%
\special{dt 0.045}%
\special{pa 5707 2815}%
\special{pa 5706 2815}%
\special{dt 0.045}%
\special{pa 5707 2815}%
\special{pa 5707 2695}%
\special{dt 0.045}%
\special{pa 5707 2695}%
\special{pa 5707 2696}%
\special{dt 0.045}%
%
\special{pn 4}%
\special{pa 5830 2753}%
\special{pa 5750 2833}%
\special{fp}%
\special{pa 5840 2803}%
\special{pa 5790 2853}%
\special{fp}%
\special{pa 5780 2743}%
\special{pa 5710 2813}%
\special{fp}%
%
\special{pn 4}%
\special{pa 5700 2823}%
\special{pa 5650 2873}%
\special{fp}%
\special{pa 5750 2833}%
\special{pa 5710 2873}%
\special{fp}%
%
\special{pn 4}%
\special{pa 5650 2753}%
\special{pa 5580 2823}%
\special{fp}%
\special{pa 5600 2743}%
\special{pa 5560 2783}%
\special{fp}%
%
\special{pn 4}%
\special{pa 5710 2753}%
\special{pa 5670 2793}%
\special{fp}%
%
\special{pn 4}%
\special{pa 5700 2703}%
\special{pa 5660 2743}%
\special{fp}%
%
\special{pn 4}%
\special{pa 5750 2713}%
\special{pa 5720 2743}%
\special{fp}%
\put(18.9000,-7.4000){\makebox(0,0)[lb]{$L_1$}}%
\put(14.6000,-12.0000){\makebox(0,0)[lb]{$L_2$}}%
\put(13.4000,-17.3000){\makebox(0,0)[lb]{$V_1$}}%
\put(20.1000,-25.8000){\makebox(0,0)[lb]{$S_1$}}%
\put(19.4000,-29.0000){\makebox(0,0)[lb]{$S_2$}}%
\put(32.1000,-6.7000){\makebox(0,0)[lb]{$V_2$}}%
\put(31.7000,-16.0000){\makebox(0,0)[lb]{$S_3$}}%
\put(26.2000,-26.2000){\makebox(0,0)[lb]{$V_3$}}%
%
\special{pn 4}%
\special{pa 4440 560}%
\special{pa 4340 660}%
\special{fp}%
\special{pa 4420 520}%
\special{pa 4330 610}%
\special{fp}%
\special{pa 4390 490}%
\special{pa 4330 550}%
\special{fp}%
\special{pa 4370 450}%
\special{pa 4330 490}%
\special{fp}%
\special{pa 4470 590}%
\special{pa 4400 660}%
\special{fp}%
\special{pa 4490 630}%
\special{pa 4460 660}%
\special{fp}%
%
\special{pn 4}%
\special{pa 4520 660}%
\special{pa 4430 750}%
\special{fp}%
\special{pa 4510 730}%
\special{pa 4470 770}%
\special{fp}%
\special{pa 4460 660}%
\special{pa 4390 730}%
\special{fp}%
\special{pa 4400 660}%
\special{pa 4340 720}%
\special{fp}%
\put(44.6000,-6.1000){\makebox(0,0)[lb]{$S_4$}}%
\put(39.2000,-16.2000){\makebox(0,0)[lb]{$V_4$}}%
\put(39.8000,-29.1000){\makebox(0,0)[lb]{$S_5$}}%
\put(58.1000,-9.0000){\makebox(0,0)[lb]{$V_5$}}%
\put(52.6000,-19.5000){\makebox(0,0)[lb]{$S_6$}}%
\put(58.3000,-29.0000){\makebox(0,0)[lb]{$V_6$}}%
%
\special{pn 20}%
\special{pa 5640 2280}%
\special{pa 5640 2088}%
\special{fp}%
\special{sh 1}%
\special{pa 5640 2088}%
\special{pa 5620 2155}%
\special{pa 5640 2141}%
\special{pa 5660 2155}%
\special{pa 5640 2088}%
\special{fp}%
\end{picture}%
\label{fig:F11}
\caption{This figure denotes the steps which are needed to obtain the 4-fold $\mathcal S$. The first figure denotes the boundary divisor ${\Bbb P}^3$ in ${\Bbb P}^4$. Up arrow denotes blowing-up, and down arrow denotes blowing-down. Each step is explained in the below summary.}
\end{figure}
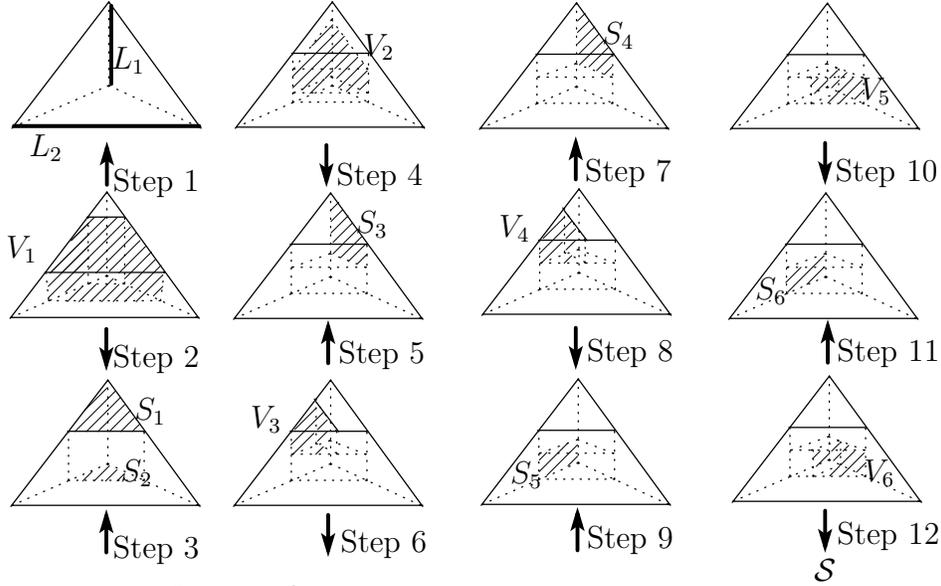

Let us summarize the steps which are needed to obtain the 4-fold $\mathcal S$.
\begin{enumerate}
\item Blow up along two curves $L_1 \cong {\Bbb P}^1$ and $L_2 \cong {\Bbb P}^1$.
\item Blow down the 3-fold $V_1 \cong {\Bbb P}^1 \times {\Bbb P}^1 \times {\Bbb P}^1$.
\item Blow up along two surfaces $S_1 \cong {\Bbb P}^2$ and $S_2 \cong {\Bbb P}^2$.
\item Blow down the 3-fold $V_2 \cong {\Bbb P}^2 \times {\Bbb P}^1$.
\item Blow up along the surface $S_3 \cong {\Bbb P}^1 \times {\Bbb P}^1$.
\item Blow down the 3-fold $V_3 \cong {\Bbb P}^2 \times {\Bbb P}^1$.
\item Blow up along the surface $S_4 \cong {\Bbb P}^1 \times {\Bbb P}^1$.
\item Blow down the 3-fold $V_4 \cong {\Bbb P}^2 \times {\Bbb P}^1$.
\item Blow up along the surface $S_5 \cong {\Bbb P}^1 \times {\Bbb P}^1$.
\item Blow down the 3-fold $V_5 \cong {\Bbb P}^2 \times {\Bbb P}^1$.
\item Blow up along the surface $S_6 \cong {\Bbb P}^1 \times {\Bbb P}^1$.
\item Blow down the 3-fold $V_6 \cong {\Bbb P}^2 \times {\Bbb P}^1$.
\end{enumerate}

It is easy to see that this rational vector field $\tilde v$ satisfies the condition:
\begin{equation}
\tilde v \in H^0({\mathcal S},\Theta_{\mathcal S}(-\log{\mathcal H})({\mathcal H})).
\end{equation}

The following lemma shows that this rational vector field $\tilde v$ has five accessible singular loci on the boundary divisor $\mathcal H \times \{ $t$ \} \subset {\mathcal S} \times \{ $t$ \}$ for each $t \in B$.

\begin{figure}[h]
\unitlength 0.1in
\begin{picture}(49.60,23.60)(19.80,-24.90)
%
\special{pn 8}%
\special{pa 4475 400}%
\special{pa 2170 2191}%
\special{fp}%
%
\special{pn 8}%
\special{pa 2192 2177}%
\special{pa 6661 2177}%
\special{fp}%
%
\special{pn 8}%
\special{pa 4475 400}%
\special{pa 6661 2191}%
\special{fp}%
%
\special{pn 8}%
\special{pa 4475 400}%
\special{pa 4475 1599}%
\special{dt 0.045}%
\special{pa 4475 1599}%
\special{pa 4475 1598}%
\special{dt 0.045}%
\special{pa 2192 2177}%
\special{pa 4475 1584}%
\special{dt 0.045}%
\special{pa 4475 1584}%
\special{pa 4474 1584}%
\special{dt 0.045}%
\special{pa 4475 1584}%
\special{pa 6616 2177}%
\special{dt 0.045}%
\special{pa 6616 2177}%
\special{pa 6615 2177}%
\special{dt 0.045}%
%
\special{pn 8}%
\special{pa 3513 1125}%
\special{pa 5341 1125}%
\special{fp}%
%
\special{pn 8}%
\special{pa 3540 1125}%
\special{pa 3540 1821}%
\special{dt 0.045}%
\special{pa 3540 1821}%
\special{pa 3540 1820}%
\special{dt 0.045}%
\special{pa 3540 1821}%
\special{pa 5336 1821}%
\special{dt 0.045}%
\special{pa 5336 1821}%
\special{pa 5335 1821}%
\special{dt 0.045}%
\special{pa 5336 1821}%
\special{pa 5336 1140}%
\special{dt 0.045}%
\special{pa 5336 1140}%
\special{pa 5336 1141}%
\special{dt 0.045}%
%
\special{pn 8}%
\special{pa 4456 1259}%
\special{pa 3558 1451}%
\special{dt 0.045}%
\special{pa 3558 1451}%
\special{pa 3559 1451}%
\special{dt 0.045}%
\special{pa 3558 1451}%
\special{pa 5314 1451}%
\special{dt 0.045}%
\special{pa 5314 1451}%
\special{pa 5313 1451}%
\special{dt 0.045}%
\special{pa 5314 1451}%
\special{pa 4456 1259}%
\special{dt 0.045}%
\special{pa 4456 1259}%
\special{pa 4457 1259}%
\special{dt 0.045}%
%
\special{pn 20}%
\special{pa 2674 1786}%
\special{pa 3538 1561}%
\special{fp}%
%
\special{pn 20}%
\special{pa 2395 2002}%
\special{pa 3538 1714}%
\special{fp}%
%
\special{pn 20}%
\special{pa 4276 1453}%
\special{pa 4474 2164}%
\special{fp}%
%
\special{pn 20}%
\special{pa 6643 2173}%
\special{pa 5356 1129}%
\special{fp}%
%
\special{pn 20}%
\special{pa 6355 2092}%
\special{pa 5329 1300}%
\special{fp}%
\put(22.6900,-19.6600){\makebox(0,0)[lb]{$C_0$}}%
\put(24.3100,-17.4100){\makebox(0,0)[lb]{$C_1$}}%
\put(42.8500,-23.2600){\makebox(0,0)[lb]{$C_2$}}%
\put(59.5900,-21.3700){\makebox(0,0)[lb]{$C_3$}}%
\put(65.1000,-20.7000){\makebox(0,0)[lb]{$C_4$}}%
\put(21.6000,-24.3000){\makebox(0,0)[lb]{$\frac{1}{Y_3}$}}%
\put(45.5000,-4.1000){\makebox(0,0)[lb]{$\frac{1}{X_1}$}}%
\put(45.9000,-16.9000){\makebox(0,0)[lb]{$\frac{1}{Z_2}$}}%
\put(65.7000,-24.2000){\makebox(0,0)[lb]{$\frac{1}{W_4}$}}%
%
\special{pn 20}%
\special{pa 4470 400}%
\special{pa 4640 130}%
\special{fp}%
\special{sh 1}%
\special{pa 4640 130}%
\special{pa 4588 176}%
\special{pa 4612 175}%
\special{pa 4621 197}%
\special{pa 4640 130}%
\special{fp}%
%
\special{pn 20}%
\special{pa 2200 2180}%
\special{pa 1980 2480}%
\special{fp}%
\special{sh 1}%
\special{pa 1980 2480}%
\special{pa 2036 2438}%
\special{pa 2012 2437}%
\special{pa 2003 2414}%
\special{pa 1980 2480}%
\special{fp}%
%
\special{pn 20}%
\special{pa 6650 2180}%
\special{pa 6940 2490}%
\special{fp}%
\special{sh 1}%
\special{pa 6940 2490}%
\special{pa 6909 2428}%
\special{pa 6904 2451}%
\special{pa 6880 2455}%
\special{pa 6940 2490}%
\special{fp}%
%
\special{pn 20}%
\special{pa 4480 1570}%
\special{pa 4700 1470}%
\special{fp}%
\special{sh 1}%
\special{pa 4700 1470}%
\special{pa 4631 1479}%
\special{pa 4651 1492}%
\special{pa 4648 1516}%
\special{pa 4700 1470}%
\special{fp}%
\end{picture}%
\label{fig:F2}
\caption{This figure denotes the boundary divisor $\mathcal H$ of $\mathcal S$. This divisor is covered by eight affine spaces $U_3 \cup U_4 \cup U_6 \cup U_7 \cup \dots \cup U_{11}$. The bold lines $C_i \ (i=0,1,\ldots,4)$ in $\mathcal H$ denote the accessible singular loci of the system \eqref{1} (see Lemma 5.2).}
\end{figure}
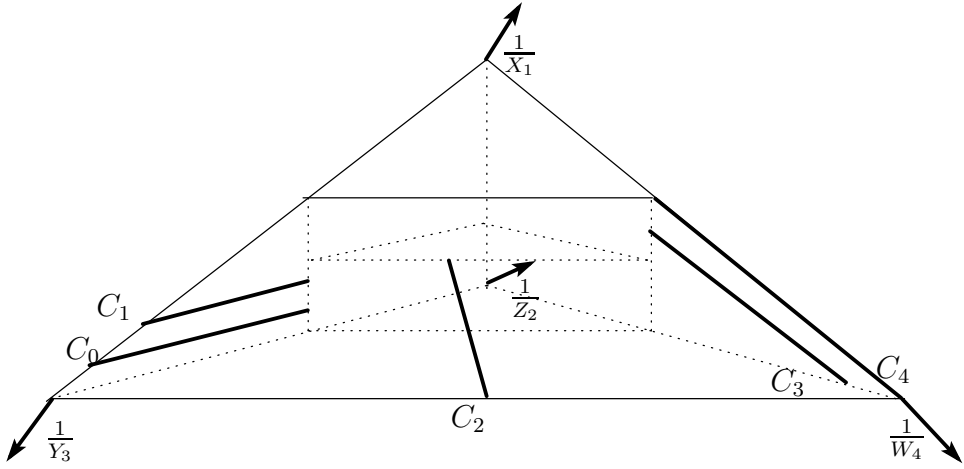

\begin{lmm}\label{9.1}

The rational vector field $\tilde v$ has the following accessible singular loci{\rm : \rm}
\begin{equation}
  \left\{
  \begin{aligned}
   C_0 &=\{(X_3,Y_3,Z_3,W_3)|X_3=t,Y_3=W_3=0\},\\
   C_1 &=\{(X_3,Y_3,Z_3,W_3)|X_3=\eta,Y_3=W_3=0\},\\
   C_2 &=\{(X_3,Y_3,Z_3,W_3)|X_3=Z_3,Y_3=0,W_3=-1\},\\
   C_3 &=\{(X_4,Y_4,Z_4,W_4)|Y_4=W_4=0,Z_4=1\},\\
   C_4 &=\{(X_4,Y_4,Z_4,W_4)|Y_4=Z_4=W_4=0\}.\\
   \end{aligned}
  \right. 
\end{equation}

\end{lmm}

This lemma can be proven by a direct calculation. \qed

Next let us calculate its local index at each point of $C_i$.

\begin{center}
\begin{tabular}{|c|c|c|} \hline 
Singular locus & Singular point & Type of local index   \\ \hline 
$C_0$ & $(X_3,Y_3,Z_3,W_3)=(t,0,a,0)$   & $(2,1,0,1)$  \\ \hline 
$C_1$ & $(X_3,Y_3,Z_3,W_3)=(\eta,0,a,0)$   & $(2,1,0,1)$  \\ \hline 
$C_2$ & $(X_4,Y_4,Z_4,W_4)=(a,-1,a,0)$   & $(0,1,2,1)$  \\ \hline 
$C_3$ & $(X_4,Y_4,Z_4,W_4)=(a,0,1,0)$   & $(0,1,2,1)$  \\ \hline 
$C_4$ & $(X_4,Y_4,Z_4,W_4)=(a,0,0,0)$   & $(0,1,2,1)$  \\ \hline 
\end{tabular}
\end{center}
Here $a \in {\Bbb C}$.

\begin{exa}
Let us take the coordinate system $(x,y,z,w)$ centered at the point $(X_3,Y_3,Z_3,W_3)=(t,0,0,0)$. The system \eqref{1} is rewritten as follows:
\begin{align*}
\frac{d}{dt}\begin{pmatrix}
             x \\
             y \\
             z \\
             w 
             \end{pmatrix}&=\frac{1}{y}\{\begin{pmatrix}
             2 & 0 & 0 & 0 \\
             0 & 1 & 0 & 0 \\
             0 & 0 & 0 & 0 \\
             0 & 0 & 0 & 1
             \end{pmatrix}\begin{pmatrix}
             x \\
             y \\
             z \\
             w 
             \end{pmatrix}+\dots\}
             \end{align*}
satisfying \eqref{b}.
\end{exa}

\begin{exa}
Let us take the coordinate system $(x,y,z,w)$ centered at the point $(X_4,Y_4,Z_4,W_4)=(0,-1,0,0)$. The system \eqref{1} is rewritten as follows:
\begin{align*}
\frac{d}{dt}\begin{pmatrix}
             x \\
             y \\
             z \\
             w 
             \end{pmatrix}&=\frac{1}{w}\{\frac{\eta}{(t-1)(t-\eta)}\begin{pmatrix}
             2 & 0 & -2 & 0 \\
             -2 & 1 & 2 & 0 \\
             0 & 0 & 0 & 0 \\
             0 & 0 & 0 & 1
             \end{pmatrix}\begin{pmatrix}
             x \\
             y \\
             z \\
             w 
             \end{pmatrix}+\dots\}
             \end{align*}
satisfying \eqref{b}. To the above system, we make the linear transformation
\begin{equation*}
\begin{pmatrix}
             X \\
             Y \\
             Z \\
             W 
             \end{pmatrix}=\begin{pmatrix}
             0 & 0 & 1 & 0 \\
             0 & 0 & 0 & 1 \\
             -1 & 0 & 1 & 0 \\
             2 & 1 & -2 & 0
             \end{pmatrix}\begin{pmatrix}
             x \\
             y \\
             z \\
             w 
             \end{pmatrix}
\end{equation*}
to arrive at
\begin{equation*}
\frac{d}{dt}\begin{pmatrix}
             X \\
             Y \\
             Z \\
             W 
             \end{pmatrix}=\frac{1}{W}\{\frac{\eta}{(t-1)(t-\eta)}\begin{pmatrix}
             0 & 0 & 0 & 0 \\
             0 & 1 & 0 & 0 \\
             0 & 0 & 2 & 0 \\
             0 & 0 & 0 & 1
             \end{pmatrix}\begin{pmatrix}
             X \\
             Y \\
             Z \\
             W 
             \end{pmatrix}+\dots\}.
             \end{equation*}
\end{exa}

\begin{prp}\label{pro:ini}\label{prop}
If we resolve the accessible singular loci given in Lemma \ref{9.1} by blowing-ups, then we can obtain the canonical coordinates $r_j (j=0,1,3,5,6)$.
\end{prp}

{\it Proof \ref{prop}.} \hspace{0.2cm} By the following steps, we can resolve the accessible singular locus $C_4$.

{\bf Step 1}: We blow up along the curve $C_4${\rm : \rm}
$$
{X_4}^{(1)}=X_4 \;, \;\;\; {Y_4}^{(1)}=\frac{Y_4}{W_4} \;, \;\;\; {Z_4}^{(1)}=\frac{Z_4}{W_4} \;, \;\;\; {W_4}^{(1)}=W_4.
$$

{\bf Step 2}: We blow up along the surface $\{({X_4}^{(1)},{Y_4}^{(1)},{Z_4}^{(1)},{W_4}^{(1)})|{Z_4}^{(1)}-\alpha_6\\
={W_4}^{(1)}=0 \}${\rm : \rm}
$$
{X_4}^{(2)}={X_4}^{(1)} \;, \;\;\; {Y_4}^{(2)}={Y_4}^{(1)} \;, \;\;\; {Z_4}^{(2)}=\frac{{Z_4}^{(1)}-\alpha_6}{{W_4}^{(1)}} \;, \;\;\; {W_4}^{(2)}={W_4}^{(1)}.
$$
Thus we have resolved the accessible singular locus $C_4$. 

By choosing a new coordinate system as
\begin{equation*}
(x_6,y_6,z_6,w_6)=({X_4}^{(2)},{Y_4}^{(2)},-{Z_4}^{(2)},{W_4}^{(2)}),
\end{equation*}
we can obtain the coordinate $r_6$.

\vspace{0.5cm}
By the following steps, we can resolve the accessible singular locus $C_2$.

{\bf Step 1}: We blow up along the curve $C_2${\rm : \rm}
$$
{X_5}^{(1)}=\frac{X_3-Z_3}{Y_3} \;, \;\;\; {Y_5}^{(1)}=Y_3 \;, \;\;\; {Z_5}^{(1)}=Z_3 \;, \;\;\; {W_5}^{(1)}=\frac{W_3+1}{Y_3}.
$$

{\bf Step 2}: We blow up along the surface $\{({X_5}^{(1)},{Y_5}^{(1)},{Z_5}^{(1)},{W_5}^{(1)})| {X_5}^{(1)}-\alpha_3={Y_5}^{(1)}=0\}${\rm : \rm}
$$
{X_5}^{(2)}=\frac{{X_5}^{(1)}-\alpha_3}{{Y_5}^{(1)}} \;, \;\;\; {Y_5}^{(2)}={Y_5}^{(1)} \;, \;\;\; {Z_5}^{(2)}={Z_5}^{(1)} \;, \;\;\; {W_5}^{(2)}={W_5}^{(1)}.
$$
Thus we have resolved the accessible singular locus $C_2$. 

By choosing a new coordinate system as
$$
(x_3,y_3,z_3,w_3)=(-{X_5}^{(2)},{Y_5}^{(2)},{Z_5}^{(2)},{W_5}^{(2)}),
$$
we can obtain the coordinate $r_3$.

For the remaining accessible singular locus, the proof is similar.

Collecting all the cases, we have obtained the canonical coordinate systems $(x_j,y_j,\\
z_j,w_j) \ (j=0,1,3,5,6)$, which proves Proposition 5.5.  \qed

We remark that each coordinate system contains a three-parameter family of meromorphic solutions of \eqref{1} as the initial conditions.

\begin{figure}[h]
\unitlength 0.1in
\begin{picture}(35.93,14.33)(21.70,-18.33)
%
\special{pn 8}%
\special{pa 4014 400}%
\special{pa 2170 1833}%
\special{fp}%
%
\special{pn 8}%
\special{pa 2188 1822}%
\special{pa 5763 1822}%
\special{fp}%
%
\special{pn 8}%
\special{pa 4014 400}%
\special{pa 5763 1833}%
\special{fp}%
%
\special{pn 8}%
\special{pa 4014 400}%
\special{pa 4014 1359}%
\special{dt 0.045}%
\special{pa 4014 1359}%
\special{pa 4014 1358}%
\special{dt 0.045}%
\special{pa 2188 1822}%
\special{pa 4014 1347}%
\special{dt 0.045}%
\special{pa 4014 1347}%
\special{pa 4013 1347}%
\special{dt 0.045}%
\special{pa 4014 1347}%
\special{pa 5727 1822}%
\special{dt 0.045}%
\special{pa 5727 1822}%
\special{pa 5726 1822}%
\special{dt 0.045}%
%
\special{pn 8}%
\special{pa 3244 980}%
\special{pa 4707 980}%
\special{fp}%
%
\special{pn 8}%
\special{pa 3266 980}%
\special{pa 3266 1537}%
\special{dt 0.045}%
\special{pa 3266 1537}%
\special{pa 3266 1536}%
\special{dt 0.045}%
\special{pa 3266 1537}%
\special{pa 4703 1537}%
\special{dt 0.045}%
\special{pa 4703 1537}%
\special{pa 4702 1537}%
\special{dt 0.045}%
\special{pa 4703 1537}%
\special{pa 4703 992}%
\special{dt 0.045}%
\special{pa 4703 992}%
\special{pa 4703 993}%
\special{dt 0.045}%
%
\special{pn 8}%
\special{pa 3999 1087}%
\special{pa 3280 1241}%
\special{dt 0.045}%
\special{pa 3280 1241}%
\special{pa 3281 1241}%
\special{dt 0.045}%
\special{pa 3280 1241}%
\special{pa 4685 1241}%
\special{dt 0.045}%
\special{pa 4685 1241}%
\special{pa 4684 1241}%
\special{dt 0.045}%
\special{pa 4685 1241}%
\special{pa 3999 1087}%
\special{dt 0.045}%
\special{pa 3999 1087}%
\special{pa 4000 1087}%
\special{dt 0.045}%
%
\special{pn 20}%
\special{pa 2573 1509}%
\special{pa 3264 1329}%
\special{fp}%
\put(23.7900,-14.7300){\makebox(0,0)[lb]{$C_1$}}%
%
\special{pn 20}%
\special{pa 3266 992}%
\special{pa 3266 1240}%
\special{fp}%
\put(30.9000,-9.7600){\makebox(0,0)[lb]{$C_{\infty}$}}%
%
\special{pn 8}%
\special{pa 2980 1340}%
\special{pa 2984 1306}%
\special{pa 2990 1274}%
\special{pa 2999 1244}%
\special{pa 3012 1217}%
\special{pa 3031 1194}%
\special{pa 3055 1174}%
\special{pa 3082 1158}%
\special{pa 3113 1143}%
\special{pa 3145 1130}%
\special{pa 3170 1120}%
\special{sp}%
%
\special{pn 8}%
\special{pa 3130 1140}%
\special{pa 3220 1080}%
\special{dt 0.045}%
\special{sh 1}%
\special{pa 3220 1080}%
\special{pa 3153 1100}%
\special{pa 3176 1110}%
\special{pa 3176 1134}%
\special{pa 3220 1080}%
\special{fp}%
\end{picture}%
\label{fig:F3}
\caption{}
\end{figure}

The difference between $r_i$ and $r_i'$ is only the case of $i=1$. The relation between $r_1$ and $r_1'$ can be explained by the one for the accessible singularities $C_1$ and $C_{\infty}$ given by
\begin{align}
\begin{split}
   C_1 =&\{(X_6,Y_6,Z_6,W_6)|X_6=\frac{1}{\eta},Y_3=W_3=0\}\\
    &\cup \{(X_8,Y_8,Z_8,W_8)|X_8=\frac{1}{\eta},Y_8=W_8=0\},\\
    C_{\infty} =&\{(X_6,Y_6,Z_6,W_6)|X_6=Y_3=W_3=0\}\\
    &\cup \{(X_8,Y_8,Z_8,W_8)|X_8=Y_8=W_8=0\}.
    \end{split}
\end{align}
As $\eta \rightarrow \infty$, $C_1$ tends to $C_{\infty}$. The resolution of $C_{\infty}$ is the same way given in Proof 5.5.

\begin{prp}
After a series of explicit blowing-ups given in Proposition \ref{prop}, we obtain the smooth projective 4-fold $\tilde{\mathcal S}$ and a morphism $\varphi:\tilde{\mathcal S} \rightarrow \mathcal S$. Its canonical divisor $K_{\tilde{\mathcal S}}$ of $\tilde{\mathcal S}$ is given by
\begin{align}
\begin{split}
K_{\tilde{\mathcal S}}&=-3{\tilde{\mathcal H}} - \sum_{i=0}^{4} {\mathcal E}_i,
\end{split}
\end{align}
where the symbol $\tilde{\mathcal H}$ denotes the proper transform of $\mathcal H$ by $\varphi$ and  ${\mathcal E}_i$ denote the exceptional divisors obtained by Step 1 (see Proof of Proposition \ref{prop}).
\end{prp}

{\it Acknowledgement.} The author would like to thank the referee, H. Kawamuko, M. Murata and M. Noumi for useful comments.

\end{document}